\documentclass[11pt]{amsart} 
\usepackage{amsmath,amssymb,mathrsfs}
\newtheorem{theorem}{Theorem}[section]
\newtheorem{proposition}[theorem]{Proposition}
\newtheorem{corollary}[theorem]{Corollary}
\newtheorem{lemma}[theorem]{Lemma}
\theoremstyle{definition}
\newtheorem{definition}[theorem]{Definition}
\newtheorem{assumption}[theorem]{Assumption}

\begin{document}

\title{Gluing Eguchi-Hanson metrics and a question of Page}
\author{Simon Brendle and Nikolaos Kapouleas}
\address{Department of Mathematics \\ Stanford University \\ Stanford, CA 94305}
\address{Department of Mathematics \\ Brown University \\ Providence, RI 02912}
\thanks{The first author was supported in part by the National Science Foundation under grants DMS-0905628 and DMS-1201924. The second author was supported in part by the National Science Foundation under grant DMS-1105371.}

\begin{abstract}
In 1978, Gibbons-Pope and Page proposed a physical picture for the Ricci flat K\"ahler metrics on the K3 surface based on a gluing construction. In this construction, one starts from a flat torus with $16$ orbifold points, and resolves the orbifold singularities by gluing in $16$ small Eguchi-Hanson manifolds which all have the same orientation. This construction was carried out rigorously by Topiwala, LeBrun-Singer, and Donaldson.

In 1981, Page asked whether the above construction can be modified by reversing the orientations of some of the Eguchi-Hanson manifolds. This is a subtle question: if successful, this construction would produce Einstein metrics which are neither K\"ahler nor self-dual. 

In this paper, we focus on a configuration of maximal symmetry involving $8$ small Eguchi-Hanson manifolds of each orientation which are arranged according to a chessboard pattern. By analyzing the interactions between Eguchi-Hanson manifolds with opposite orientation, we identify a non-vanishing obstruction to the gluing problem, thereby destroying any hope of producing a metric of zero Ricci curvature in this way. Using this obstruction, we are able to understand the dynamics of such metrics under Ricci flow as long as the Eguchi-Hanson manifolds remain small. In particular, for the configuration described above, we obtain an ancient solution to the Ricci flow with the property that the maximum of the Riemann curvature tensor blows up at a rate of $(-t)^{\frac{1}{2}}$, while the maximum of the Ricci curvature converges to $0$.

\end{abstract}
\maketitle

\section{Introduction}

Gluing techniques are a central tool for constructing solutions of nonlinear partial differential equations. In a first step, one constructs an approximate solution of the given partial differential equation by gluing together. In the next step, one tries to deform this approximate solution to an exact solution using the implicit function theorem. This method has been used to construct solutions of many nonlinear elliptic equations arising in differential geometry (see e.g. \cite{Anderson2}, \cite{Bamler}, \cite{Biquard-Minerbe}, \cite{Brendle0}, \cite{Gursky-Viaclovsky}, \cite{Joyce}, \cite{Kapouleas1}, \cite{Kapouleas2}, \cite{Kapouleas3}, \cite{Kapouleas-survey}, \cite{Kapouleas4}, \cite{Mazzeo-Pacard1}, \cite{Mazzeo-Pacard2}, \cite{Schoen}, \cite{Taubes}, \cite{Taubes-survey}). 

In this paper, we consider a gluing problem for the Einstein equations and their parabolic analogue, the Ricci flow (cf. \cite{Hamilton1}, \cite{Hamilton2}, \cite{Hamilton-survey}; see also \cite{Brendle-book}). Our starting point is the torus $\mathbb{R}^4 / (2\mathbb{Z})^4$ equipped with the flat metric. After identifying each point on $\mathbb{R}^4 / (2\mathbb{Z})^4$ with its image under antipodal reflection, we obtain a flat orbifold with $16$ singular points. We may desingularize this orbifold by gluing in an Eguchi-Hanson manifold at each of the $16$ singular points. Recall that the Eguchi-Hanson manifold is a smooth, Ricci flat K\"ahler manifold. The Eguchi-Hanson metric on the complement of a two-sphere can be written in coordinates as 
\[g_{\text{\rm eh},\varepsilon} = \frac{r^2}{(\varepsilon^4+r^4)^{\frac{1}{2}}} \, (dr \otimes dr+r^2 \, \alpha_1 \otimes \alpha_1) + (\varepsilon^4+r^4)^{\frac{1}{2}} \, (\alpha_2 \otimes \alpha_2+\alpha_3 \otimes \alpha_3).\] 
Here, $r$ denotes the radial coordinate, $\alpha_1,\alpha_2,\alpha_3$ is a suitable set of one-forms on $S^3$, 
and $\varepsilon$ is a scaling parameter. 
An important point is that the Eguchi-Hanson manifold is asymptotic to $\mathbb{R}^4/\mathbb{Z}_2$ near infinity. 
Thus, inserting an Eguchi-Hanson at each orbifold point will result in a smooth manifold. 
Since the torus and the Eguchi-Hanson are all Ricci flat, the resulting manifold admits a metric which is almost Ricci flat. 

The key issue is whether this approximate solution of the Einstein equation can be deformed to an exact one. It turns out that the choice of orientation of the $16$ Eguchi-Hanson manifolds plays a crucial role. In the special case when the $16$ Eguchi-Hanson manifolds all have the same orientation, Gibbons and Pope \cite{Gibbons-Pope} and Page \cite{Page1} suggested that this gluing construction should recover the Ricci flat K\"ahler metric on the K3 surface. Rigorous proofs of this fact were given by Topiwala \cite{Topiwala1},\cite{Topiwala2}, LeBrun and Singer \cite{LeBrun-Singer}, and Donaldson \cite{Donaldson}. An important point here is that one can work within the class of K\"ahler manifolds. Thus, the problem of finding a Ricci flat metric can be reduced to the solvability of a Monge-Amp\`ere-type equation. This construction also shares some common features with Joyce's construction of manifolds with exceptional holonomy (cf. \cite{Joyce}).  

In the following, we consider a gluing construction involving Eguchi-Hanson manifolds with different orientations. Such a construction was first proposed by Page \cite{Page2} in 1981. One major difficulty in this case is that the resulting metric will not be K\"ahler, and it is necessary to work with the full Einstein equations. For simplicity, we consider a particularly symmetric configuration involving $8$ Eguchi-Hanson manifolds with positive orientation and $8$ Eguchi-Hanson manifolds of negative orientation, where the orientations are assigned according to a checkboard pattern. It turns out that, modulo symmetries, the approximate kernel of the linearized operator is two-dimensional. From a geometric point of view, the approximate kernel reflects the freedom to change the scaling parameter $\varepsilon$ of the Eguchi-Hanson metrics and the size of the flat metric on the torus. The crucial issue then is to analyze the projection of the Ricci tensor to the approximate kernel. It turns out that this projection is non-zero, due to interactions between Eguchi-Hanson metrics with opposite orientations. These interactions have the effect of increasing the size of the Eguchi-Hanson metrics when the metric is evolved by the Ricci flow. In fact, a formal calculation suggests that, under the Ricci flow, the scale parameter $\varepsilon$ should evolve according to the equation $\frac{d}{dt} \varepsilon = 8\omega \, \varepsilon^5 + o(\varepsilon^5)$, where $\omega$ is a positive constant defined in Proposition \ref{flux.integral}. This is reminiscent of the work of Struwe \cite{Struwe} on the heat flow for the Nirenberg problem.

In light of this obstruction, we cannot expect to deform the metric to one of zero Ricci curvature. Instead, we show that this setup leads to a non-trivial ancient solution to the Ricci flow:

\begin{theorem}
\label{main.thm}
There exists a compact ancient solution to the Ricci flow in dimension $4$ with the following property.  For $-t$ sufficiently large, the manifold can be viewed as a desingularization of a flat torus with $16$ orbifold points. More precisely, we divide the $16$ orbifold points into two classes according to a checkboard pattern. Near $8$ orbifold points, the metric is a small perturbation of a positively-oriented near the remaining $8$ orbifold points, the metric is a small perturbation of a negatively-oriented Eguchi-Hanson metric. As $t \to -\infty$, the size of the Eguchi-Hanson instantons shrinks to zero, and we have $\sup |\text{\rm Rm}_{g(t)}|_{g(t)} = (c+o(1)) \, (-t)^{\frac{1}{2}}$, where $c$ is a positive constant. Finally, the Ricci curvature of $g(t)$ satisfies $\sup |\text{\rm Ric}_{g(t)}|_{g(t)} = O((-t)^{-\frac{1}{2}+\kappa})$ as $t \to -\infty$, where $\kappa>0$ can be chosen arbitrarily small.
\end{theorem}

Ancient solutions play a crucial role as singularity models for the Ricci flow; see e.g. \cite{Hamilton-survey}, \cite{Perelman1}, and \cite{Perelman2}. It is an interesting question to classify ancient solutions and Ricci solitons. We refer to \cite{Brendle1}, \cite{Brendle2}, \cite{Brendle-Huisken-Sinestrari}, \cite{Daskalopoulos-Hamilton-Sesum} for some recent progress in this direction. 

An interesting open question is what happens to our ancient solution beyond the range of sufficiently large $-t$.
Another interesting observation is that if we consider less symmetric configurations we can construct metrics by gluing $16$ 
small Eguchi-Hanson metrics on the orbifold as before, which under the Ricci flow evolve so that one of the Eguchi-Hanson metrics 
becomes extinct by shrinking to zero size, while the rest stay close to fixed nonzero sizes. 
This construction will appear elsewhere. 

Theorem \ref{main.thm} is inspired in part by the remarkable recent work of Biquard. 
In \cite{Biquard}, Biquard glued an Eguchi-Hanson manifold to a given Einstein orbifold. 
The resulting metric is an approximate solution of the Einstein equation. 
Biquard found an obstruction involving the curvature of the background orbifold which, in general, prevents one from deforming this metric to an exact solution of the Einstein equation. This result can be viewed as converse of the compactness results of Anderson \cite{Anderson1} and Bando, Kasue, and Nakajima \cite{Bando-Kasue-Nakajima} (see also the survey paper \cite{Cheeger}). We note that Biquard's obstruction vanishes in our situation, as the background orbifold is a flat torus. 
For that reason, we need to perform a more precise calculation which takes into account the interactions between different Eguchi-Hanson manifolds.
Note that the construction and proof are modelled after a gluing construction for an elliptic problem: 
We first construct an approximate ancient solution and then we prove that it can be perturbed to an exact one 
by carefully estimating the error terms and applying the Schauder fixed point theorem.

Theorem \ref{main.thm} also shares some common features with the beautiful work of Daskalopoulos, del~Pino, and \v Se\v sum \cite{Daskalopoulos-del-Pino-Sesum} on the Yamabe flow. The main result in \cite{Daskalopoulos-del-Pino-Sesum} asserts that there exists a non-trivial ancient solution to the Yamabe flow which is conformally flat and looks like two spheres joined by a small neck when $-t$ is large. 

We also note that Takahashi \cite{Takahashi} has recently constructed an ancient solution to the Ricci flat which converges smoothly to the Euclidean Schwarzschild metric as $t \to -\infty$. In particular, this solution has uniformly bounded curvature as $t \to -\infty$, whereas the solution in Theorem \ref{main.thm} has unbounded curvature.

The authors are very grateful to Professors Olivier Biquard and Andr\'as Vasy for helpful discussions. 

\section{Basic properties of the Eguchi-Hanson metric}

In this section, we review the definition of the Eguchi-Hanson metric (cf. \cite{Biquard}, \cite{Eguchi-Hanson}). As in \cite{Biquard}, we define one-forms $\alpha_1,\alpha_2,\alpha_3$ on $\mathbb{R}^4 \setminus \{0\}$ by 
\begin{align*} 
\alpha_1 &= \frac{1}{r^2} \, (x_1 \, dx_2 - x_2 \, dx_1 + x_3 \, dx_4 - x_4 \, dx_3), \\ 
\alpha_2 &= \frac{1}{r^2} \, (x_1 \, dx_3 - x_3 \, dx_1 + x_4 \, dx_2 - x_2 \, dx_4), \\ 
\alpha_3 &= \frac{1}{r^2} \, (x_1 \, dx_4 - x_4 \, dx_1 + x_2 \, dx_3 - x_3 \, dx_2).
\end{align*}
The Eguchi-Hanson metric with parameter $\varepsilon>0$ is defined by 
\[g_{\text{\rm eh},\varepsilon} = \frac{r^2}{(\varepsilon^4+r^4)^{\frac{1}{2}}} \, (dr \otimes dr+r^2 \, \alpha_1 \otimes \alpha_1) + (\varepsilon^4+r^4)^{\frac{1}{2}} \, (\alpha_2 \otimes \alpha_2+\alpha_3 \otimes \alpha_3).\] 
This defines a Ricci flat metric on $\mathbb{R}^4 \setminus \{0\}$ which is invariant under antipodal reflection. The induced metric on $(\mathbb{R}^4 \setminus \{0\}) / \mathbb{Z}_2$ admits a smooth compactification, where the origin is replaced by a two-dimensional sphere. The parameter $\varepsilon$ serves as a scaling parameter; that is, different choices of $\varepsilon$ result in metrics which are isometric up to scaling.

Let $\hat{g}_{\text{\rm eh},\varepsilon}$ denote the pull-back of $g_{\text{\rm eh},\varepsilon}$ under the map $(x_1,x_2,x_3,x_4) \mapsto (-x_1,x_2,x_3,x_4)$. Clearly, $\hat{g}_{\text{\rm eh},\varepsilon}$ is a Ricci flat metric. Near infinity, we have the asymptotic expansions 
\[g_{\text{\rm eh},\varepsilon} = g_{\text{\rm eucl}} + \frac{1}{2} \, \varepsilon^4 \, T + O(\varepsilon^8 \, r^{-8})\] 
and 
\[\hat{g}_{\text{\rm eh},\varepsilon} = g_{\text{\rm eucl}} + \frac{1}{2} \, \varepsilon^4 \, \hat{T} + O(\varepsilon^8 \, r^{-8}),\] 
where 
\begin{align*} 
T &= -r^{-4} \, (dr \otimes dr + r^2 \, \alpha_1 \otimes \alpha_1 - r^2 \, \alpha_2 \otimes \alpha_2 - r^2 \, \alpha_3 \otimes \alpha_3) \\ 
&= -r^{-6} \, \big ( (x_1^2+x_2^2-x_3^2-x_4^2) \, (dx_1 \otimes dx_1 + dx_2 \otimes dx_2 - dx_3 \otimes dx_3 - dx_4 \otimes dx_4) \\ 
&\hspace{20mm} + 2 \, (x_1x_3+x_2x_4) \, (dx_1 \otimes dx_3 + dx_3 \otimes dx_1 + dx_2 \otimes dx_4 + dx_4 \otimes dx_2) \\ 
&\hspace{20mm} + 2 \, (x_1x_4-x_2x_3) \, (dx_1 \otimes dx_4 + dx_4 \otimes dx_1 - dx_2 \otimes dx_3 - dx_3 \otimes dx_2) \big ) 
\end{align*} 
and 
\begin{align*} 
\hat{T} &= -r^{-6} \, \big ( (x_1^2+x_2^2-x_3^2-x_4^2) \, (dx_1 \otimes dx_1 + dx_2 \otimes dx_2 - dx_3 \otimes dx_3 - dx_4 \otimes dx_4) \\ 
&\hspace{20mm} + 2 \, (x_1x_3-x_2x_4) \, (dx_1 \otimes dx_3 + dx_3 \otimes dx_1 - dx_2 \otimes dx_4 - dx_4 \otimes dx_2) \\ 
&\hspace{20mm} + 2 \, (x_1x_4+x_2x_3) \, (dx_1 \otimes dx_4 + dx_4 \otimes dx_1 + dx_2 \otimes dx_3 + dx_3 \otimes dx_2) \big ). 
\end{align*} 
We note that the metrics $g_{\text{\rm eh},\varepsilon}$ and $\hat{g}_{\text{\rm eh},\varepsilon}$ and the tensors $T$ and $\hat{T}$ are all invariant under the maps 
\begin{align*} 
&(x_1,x_2,x_3,x_4) \mapsto (x_2,-x_1,x_3,x_4), \\ 
&(x_1,x_2,x_3,x_4) \mapsto (x_1,x_2,x_4,-x_3), \\
&(x_1,x_2,x_3,x_4) \mapsto (x_3,x_4,x_1,x_2), \\ 
&(x_1,x_2,x_3,x_4) \mapsto (-x_3,x_4,-x_1,x_2). 
\end{align*}

In the remainder of this section, we review some known results, due to Biquard \cite{Biquard} and Page \cite{Page1}, concerning the Lichnerowicz Laplacian for the Eguchi-Hanson manifold. Recall that the Lichnerowicz Laplacian is defined by 
\[\Delta_L h_{ik} = \Delta h_{ik} + 2 R_{ijkl} h^{jl} - \text{\rm Ric}_i^l h_{kl} - \text{\rm Ric}_k^l h_{il}.\] 
To analyze the kernel of $\Delta_{L,g_{\text{\rm eh}}}$, we consider the vector fields 
\begin{align*} 
V_1 &= x_1 \, \frac{\partial}{\partial x_2} - x_2 \, \frac{\partial}{\partial x_1} + x_3 \, \frac{\partial}{\partial x_4} - x_4 \, \frac{\partial}{\partial x_3}, \\
V_2 &= x_1 \, \frac{\partial}{\partial x_3} - x_3 \, \frac{\partial}{\partial x_1} + x_4 \, \frac{\partial}{\partial x_2} - x_2 \, \frac{\partial}{\partial x_4}, \\
V_3 &= x_1 \, \frac{\partial}{\partial x_4} - x_4 \, \frac{\partial}{\partial x_1} + x_2 \, 
\frac{\partial}{\partial x_3} - x_3 \, \frac{\partial}{\partial x_2}. 
\end{align*}
Note that the frame $r \, \frac{\partial}{\partial r},V_1,V_2,V_3$ is dual to the co-frame $\frac{1}{r} \, dr,\alpha_1,\alpha_2,\alpha_3$. Moreover, $[V_1,V_2] = -2 \, V_3$, $[V_2,V_3] = -2 \, V_1$, and $[V_3,V_1] = -2 \, V_2$. This implies $\mathscr{L}_{V_1} \alpha_2 = -\mathscr{L}_{V_2} \alpha_1 = -2 \, \alpha_3$, $\mathscr{L}_{V_2} \alpha_3 = -\mathscr{L}_{V_3} \alpha_2 = -2 \, \alpha_1$, and $\mathscr{L}_{V_3} \alpha_1 = -\mathscr{L}_{V_1} \alpha_3 = -2 \, \alpha_2$. We next define 

\begin{align*} 
o_{1,\varepsilon} &= g_{\text{\rm eh},\varepsilon} - \frac{1}{2} \, \mathscr{L}_{r \frac{\partial}{\partial r}} g_{\text{\rm eh},\varepsilon}, \\ 
o_{2,\varepsilon} &= \frac{1}{2} \, \mathscr{L}_{\frac{r^2 V_2}{\sqrt{\varepsilon^4+r^4}}} g_{\text{\rm eh},\varepsilon}, \\ 
o_{3,\varepsilon} &= \frac{1}{2} \, \mathscr{L}_{\frac{r^2 V_3}{\sqrt{\varepsilon^4+r^4}}} g_{\text{\rm eh},\varepsilon}.
\end{align*} 
Equivalently, we may write 
\begin{align*} 
o_{1,\varepsilon} &= -\frac{\varepsilon^4 \, r^2}{(\varepsilon^4+r^4)^{\frac{3}{2}}} \, (dr \otimes dr + r^2 \, \alpha_1 \otimes \alpha_1) + \frac{\varepsilon^4}{(\varepsilon^4+r^4)^{\frac{1}{2}}} \, (\alpha_2 \otimes \alpha_2 + \alpha_3 \otimes \alpha_3), \\ 
o_{2,\varepsilon} &= \frac{\varepsilon^4}{\varepsilon^4+r^4} \, (r \, dr \otimes \alpha_2 + r \, \alpha_2 \otimes dr - r^2 \, \alpha_1 \otimes \alpha_3 - r^2 \, \alpha_3 \otimes \alpha_1), \\ 
o_{3,\varepsilon} &= \frac{\varepsilon^4}{\varepsilon^4+r^4} \, (r \, dr \otimes \alpha_3 + r \, \alpha_3 \otimes dr + r^2 \, \alpha_1 \otimes \alpha_2 + r^2 \, \alpha_2 \otimes \alpha_1). 
\end{align*}
Finally, the tensor $o_{1,\varepsilon}$ can be rewritten as 
\[o_{1,\varepsilon} = \frac{1}{2} \, \varepsilon \, \frac{\partial}{\partial \varepsilon} g_{\text{\rm eh},\varepsilon} = \varepsilon^4 \, T + O(\varepsilon^8 \, r^{-8}).\] 
The main properties of $o_{1,\varepsilon},o_{2,\varepsilon},o_{3,\varepsilon}$ are summarized in the following proposition.

\begin{proposition}[O.~Biquard \cite{Biquard}; D.~Page \cite{Page1}]
\label{o_i}
For each $i \in \{1,2,3\}$, the tensor $o_{i,\varepsilon}$ has the following properties: 
\begin{itemize}
\item[(i)] $\text{\rm tr}_{g_{\text{\rm eh},\varepsilon}} o_{i,\varepsilon} = 0$.  
\item[(ii)] $\text{\rm div}_{g_{\text{\rm eh},\varepsilon}} o_{i,\varepsilon} = 0$. 
\item[(iii)] $\Delta_{L,g_{\text{\rm eh},\varepsilon}} o_{i,\varepsilon} = 0$. 
\item[(iv)] $\int_{\mathbb{R}^4 \setminus \{0\}} |o_{1,\varepsilon}|_{g_{\text{\rm eh},\varepsilon}}^2 \, d\text{\rm vol}_{g_{\text{\rm eh},\varepsilon}} = 2\pi^2 \, \varepsilon^4$. 
\end{itemize}
\end{proposition}

\begin{proof} 
(i) It is obvious that $o_{i,\varepsilon}$ is trace-free for each $i \in \{1,2,3\}$. 

(ii) We first consider the vector field $Y = r \, \frac{\partial}{\partial r}$. Since $o_{1,\varepsilon} = g_{\text{\rm eh},\varepsilon} - \frac{1}{2} \, \mathscr{L}_Y g_{\text{\rm eh},\varepsilon}$ is trace-free, we have $\sum_k D_k Y^k = 4$ since $o_{1,\varepsilon}$ is trace-free. Moreover, $D_k Y_l - D_l Y_k = 0$ since $Y$ is a gradient vector field. Differentiating this identity gives $\sum_k D^k (D_k Y_l - D_l Y_k) = 0$. Therefore, we obtain 
\[\sum_k D^k (D_k Y_l + D_l Y_k) = \sum_k D^k (D_k Y_l - D_l Y_k) + 2 \sum_k D_l D^k Y_k = 0.\] 
Consequently, $o_{1,\varepsilon} = g_{\text{\rm eh},\varepsilon} - \frac{1}{2} \, \mathscr{L}_Y g_{\text{\rm eh},\varepsilon}$ is divergence-free.

In the next step, we define $Z = \frac{r^2 V_2}{\sqrt{\varepsilon^4+r^4}}$. Since $o_{2,\varepsilon} = \frac{1}{2} \, \mathscr{L}_Z g_{\text{\rm eh},\varepsilon}$ is trace-free, we have $\sum_k D_k Z^k = 0$ since $o_{2,\varepsilon}$. We next observe that $d(r^2 \, \alpha_2) = 2r \, dr \wedge \alpha_2 + 2r^2 \, \alpha_3 \wedge \alpha_1$ is a closed two-form which is self-dual with respect to the metric $g_{\text{\rm eh},\varepsilon}$. Consequently, the two-form $d(r^2 \, \alpha_2)$ is divergence-free with respect to the metric $g_{\text{\rm eh},\varepsilon}$. Since $g_{\text{\rm eh},\varepsilon}(Z,\cdot) = r^2 \, \alpha_2$, we conclude that $\sum_k D^k (D_k Z_l - D_l Z_k) = 0$. Putting these facts together, we obtain 
\[\sum_k D^k (D_k Z_l + D_l Z_k) = \sum_k D^k (D_k Z_l - D_l Z_k) + 2 \sum_k D_l D^k Z_k = 0.\] 
Thus, $o_{2,\varepsilon} = \frac{1}{2} \, \mathscr{L}_Z g_{\text{\rm eh},\varepsilon}$ is divergence-free. An analogous argument shows that $o_{3,\varepsilon}$ is divergence-free.

(iii) Clearly, $o_{i,\varepsilon}$ lies in the kernel of the linearized Einstein operator for each $i \in \{1,2,3\}$. Since $o_{i,\varepsilon}$ is trace-free and divergence-free, we conclude that $\Delta_{L,g_{\text{\rm eh},\varepsilon}} o_{i,\varepsilon} = 0$ for each $i \in \{1,2,3\}$. 

(iv) Note that 
\[|o_{1,\varepsilon}|_{g_{\text{\rm eh},\varepsilon}}^2 = 4 \, \Big ( \frac{\varepsilon^4}{\varepsilon^4+r^4} \Big )^2.\] 
This implies 
\[\int_{\mathbb{R}^4 \setminus \{0\}} |o_{1,\varepsilon}|_{g_{\text{\rm eh},\varepsilon}}^2 \, d\text{\rm vol}_{g_{\text{\rm eh},\varepsilon}} = \int_0^\infty 4 \, \Big ( \frac{\varepsilon^4}{\varepsilon^4+r^4} \Big )^2 \cdot 2\pi^2 \, r^3 \, dr = 2\pi^2 \, \varepsilon^4.\] 
This completes the proof of Proposition \ref{o_i}. 
\end{proof}

\section{Attaching Eguchi-Hanson metrics to a torus with $16$ orbifold points}

\label{attaching}

In this section, we will attach Eguchi-Hanson metrics with different orientations to a torus with $16$ orbifold points. To fix notation, we denote by $\mathscr{C}$ the collection of maps 
\begin{align*} 
&(x_1,x_2,x_3,x_4) \mapsto (x_2,-x_1,x_3,x_4), \\ 
&(x_1,x_2,x_3,x_4) \mapsto (x_1,x_2,x_4,-x_3), \\
&(x_1,x_2,x_3,x_4) \mapsto (x_3,x_4,x_1,x_2), \\ 
&(x_1,x_2,x_3,x_4) \mapsto (-x_3,x_4,-x_1,x_2), \\ 
&(x_1,x_2,x_3,x_4) \mapsto (1-x_1,x_2,x_3,x_4), \\ 
&(x_1,x_2,x_3,x_4) \mapsto (x_1,1-x_2,x_3,x_4), \\ 
&(x_1,x_2,x_3,x_4) \mapsto (x_1,x_2,1-x_3,x_4), \\ 
&(x_1,x_2,x_3,x_4) \mapsto (x_1,x_2,x_3,1-x_4), \\ 
&(x_1,x_2,x_3,x_4) \mapsto (x_1+2,x_2,x_3,x_4), \\ 
&(x_1,x_2,x_3,x_4) \mapsto (x_1,x_2+2,x_3,x_4), \\ 
&(x_1,x_2,x_3,x_4) \mapsto (x_1,x_2,x_3+2,x_4), \\ 
&(x_1,x_2,x_3,x_4) \mapsto (x_1,x_2,x_3,x_4+2), 
\end{align*} 
and by $\mathscr{G}$ the group generated by $\mathscr{C}$.

For each $a \in \mathbb{Z}^4$, we denote by $\tau_a$ the translation $x \mapsto x-a$. Moreover, we put 
\[\mathbb{Z}_{\text{\rm even}}^4 = \{(a_1,a_2,a_3,a_4) \in \mathbb{Z}^4: \text{\rm $a_1+a_2+a_3+a_4$ is even}\}\] 
and 
\[\mathbb{Z}_{\text{\rm odd}}^4 = \{(a_1,a_2,a_3,a_4) \in \mathbb{Z}^4: \text{\rm $a_1+a_2+a_3+a_4$ is odd}\}.\] 

\begin{lemma} 
The limits 
\[\sum_{a \in \mathbb{Z}_{\text{\rm even}}^4} \tau_a^* T := \lim_{N \to \infty} \sum_{a \in \mathbb{Z}_{\text{\rm even}}^4, \, \max \{|a_1|,|a_2|,|a_3|,|a_4|\} < N} \tau_a^* T\] 
and 
\[\sum_{a \in \mathbb{Z}_{\text{\rm odd}}^4} \tau_a^* \hat{T} := \lim_{N \to \infty} \sum_{a \in \mathbb{Z}_{\text{\rm odd}}^4, \, \max \{|a_1|,|a_2|,|a_3|,|a_4|\} < N} \tau_a^* \hat{T}\] 
exist for $x \in \mathbb{R}^4 \setminus \mathbb{Z}^4$. Neither series converges absolutely.
\end{lemma}

\begin{proof} 
Fix a point $x \in \mathbb{R}^4 \setminus \mathbb{Z}^4$. It is easy to see that 
\[\tau_{(a_1,a_2,a_3,a_4)}^* T(x) + \tau_{(a_3,-a_4,-a_1,a_2)}^* T(x) = O(|a|^{-5})\] 
as $|a| \to \infty$. Similarly, we have 
\[\tau_{(a_1,a_2,a_3,a_4)}^* \hat{T}(x) + \tau_{(a_3,a_4,-a_1,-a_2)}^* \hat{T}(x) = O(|a|^{-5})\] 
as $|a| \to \infty$ for each point $x \in [-\frac{1}{2},\frac{1}{2}]^4 \setminus \{0\}$. From this, we deduce that 
\[\bigg | \sum_{a \in \mathbb{Z}_{\text{\rm even}}^4, \, N \leq \max \{|a_1|,|a_2|,|a_3|,|a_4|\} < 2N} \tau_a^* T(x) \bigg | \leq O(N^{-1})\] 
and 
\[\bigg | \sum_{a \in \mathbb{Z}_{\text{\rm odd}}^4, \, N \leq \max \{|a_1|,|a_2|,|a_3|,|a_4|\} < 2N} \tau_a^* \hat{T}(x) \bigg | \leq O(N^{-1})\] 
for $N<N'$. From this, the assertion follows.
\end{proof}

\begin{lemma} 
The tensor 
\[\sum_{a \in \mathbb{Z}_{\text{\rm even}}^4} \tau_a^* T + \sum_{a \in \mathbb{Z}_{\text{\rm odd}}^4} \tau_a^* \hat{T}\] 
on $\mathbb{R}^4 \setminus \mathbb{Z}^4$ is invariant under the group $\mathscr{G}$ defined above.
\end{lemma}

\begin{proof} 
Consider the partial sums 
\[S^{(N)} := \sum_{a \in \mathbb{Z}_{\text{\rm even}}^4, \, \max\{|a_1|,|a_2|,|a_3|,|a_4|\} < N} \tau_a^* T + \sum_{a \in \mathbb{Z}_{\text{\rm odd}}^4, \, \max\{|a_1|,|a_2|,|a_3|,|a_4|\} < N} \tau_a^* \hat{T}.\] 
If we fix a point $x \in \mathbb{R}^4 \setminus \mathbb{Z}^4$ and a map $\varphi \in \mathscr{C}$, then we have $\varphi^* S^{(N)} - S^{(N)} = O(N^{-1})$ at the point $x$. Hence, the limit $\lim_{N \to \infty} S^{(N)}$ is invariant under each map $\varphi \in \mathscr{C}$. 
\end{proof}

\begin{definition} 
Given two positive numbers $\varepsilon$ and $\delta$ such that $\varepsilon \ll \delta^4 \ll 1$, we define a metric $\bar{g}_{\varepsilon,\delta}$ on the cube $[-\frac{1}{2},\frac{1}{2}]^4 \setminus \{0\}$ in the following way: For $|x| \leq \frac{1}{2} \, \delta$, we define 
\[\bar{g}_{\varepsilon,\delta} = g_{\text{\rm eh},\varepsilon}.\] 
Moreover, for $|x| \geq \delta$, we put 
\[\bar{g}_{\varepsilon,\delta} = g_{\text{\rm eucl}} + \frac{1}{2} \, \varepsilon^4 \sum_{a \in \mathbb{Z}_{\text{\rm even}}^4} \tau_a^* T + \frac{1}{2} \, \varepsilon^4 \sum_{a \in \mathbb{Z}_{\text{\rm odd}}^4} \tau_a^* \hat{T}.\] 
Finally, in the intermediate region $\frac{1}{2} \, \delta \leq |x| \leq \delta$, we define 
\[\bar{g}_{\varepsilon,\delta} = (1-\chi(|x|/\delta)) \, g_{\text{\rm eh},\varepsilon} + \chi(|x|/\delta) \, \bigg ( g_{\text{\rm eucl}} + \frac{1}{2} \, \varepsilon^4 \sum_{a \in \mathbb{Z}_{\text{\rm even}}^4} \tau_a^* T + \frac{1}{2} \, \varepsilon^4 \sum_{a \in \mathbb{Z}_{\text{\rm odd}}^4} \tau_a^* T \bigg ),\] 
where $\chi$ is a cutoff function satisfying $\chi=0$ on $(-\infty,\frac{2}{3}]$ and $\chi=0$ on $[\frac{5}{6},\infty)$. 
\end{definition}

We may extend $\bar{g}_{\varepsilon,\delta}$ to a metric on $\mathbb{R}^4 \setminus \mathbb{Z}^4$ which is invariant under the group $\mathscr{G}$ defined above. Note that the resulting metric $\bar{g}_{\varepsilon,\delta}$ on $\mathbb{R}^4 \setminus \mathbb{Z}^4$ is singular at each lattice point. In fact, if $a \in \mathbb{Z}_{\text{\rm even}}^4$, then we have $\bar{g}_{\varepsilon,\delta} = \tau_a^* g_{\text{\rm eh},\varepsilon}$ in a neighborhood of $a$. Similarly, if $a \in \mathbb{Z}_{\text{\rm odd}}^4$, then we have $\bar{g}_{\varepsilon,\delta} = \tau_a^* \hat{g}_{\text{\rm eh},\varepsilon}$ in a neighborhood of $a$. Thus, if we take the quotient by translations and antipodal reflection, then the metric $\bar{g}_{\varepsilon,\delta}$ descends to a smooth metric on the quotient manifold $M$.

We next estimate the Ricci curvature of $\bar{g}_{\varepsilon,\delta}$.

\begin{proposition} 
\label{ricci.tensor}
We have 
\[\text{\rm Ric}_{\bar{g}_{\varepsilon,\delta}} = 0\] 
for $|x| \leq \frac{1}{2} \, \delta$. Moreover, we have 
\[|\text{\rm Ric}_{\bar{g}_{\varepsilon,\delta}}| \leq C \, \varepsilon^4 \, \delta^{-2}\] 
for $\frac{1}{2} \, \delta \leq |x| \leq \delta$ and 
\[|\text{\rm Ric}_{\bar{g}_{\varepsilon,\delta}}| \leq C \, \varepsilon^8 \, |x|^{-10}\] 
for $|x| \geq \delta$. Analogous estimates hold for the derivatives of $\text{\rm Ric}_{\bar{g}_{\varepsilon,\delta}}$.
\end{proposition}

\begin{proof} 
Straightforward calculation.
\end{proof}

Finally, let us denote by $\bar{o}_{1,\varepsilon,\delta}$ the trace-free part of the tensor $\frac{1}{2} \, \varepsilon \, \frac{\partial}{\partial \varepsilon} \bar{g}_{\varepsilon,\delta}$ with respect to the metric $\bar{g}_{\varepsilon,\delta}$. Clearly, $\bar{o}_{1,\varepsilon,\delta} = o_{1,\varepsilon}$ for $|x| \leq \frac{1}{2} \, \delta$, so we can think of $\bar{o}_{1,\varepsilon,\delta}$ as an extension of $o_{1,\varepsilon}$ to the manifold $M$. 

\begin{proposition} 
\label{obstruction.tensor}
We have 
\[\Delta_{L,\bar{g}_{\varepsilon,\delta}} \bar{o}_{1,\varepsilon,\delta} = 0\] 
for $|x| \leq \frac{1}{2} \, \delta$. Moreover, we have 
\[|\Delta_{L,\bar{g}_{\varepsilon,\delta}} \bar{o}_{1,\varepsilon,\delta}| \leq C \, \varepsilon^4 \, \delta^{-2}\] 
for $\frac{1}{2} \, \delta \leq |x| \leq \delta$ and 
\[|\Delta_{L,\bar{g}_{\varepsilon,\delta}} \bar{o}_{1,\varepsilon,\delta}| \leq C \, \varepsilon^8 \, |x|^{-10}\] 
for $|x| \geq \delta$.
\end{proposition} 

\begin{proof} 
Again, this follows from a straightforward calculation.
\end{proof}

\section{The projection of the Ricci tensor to the approximate kernel}

\begin{lemma}
\label{distributional.laplacian}
Fix a pair of indices $i \neq j$, and let $u$ be a smooth harmonic function which is defined on the Euclidean ball $\{|x| \leq \delta\}$. Then  
\[\int_{\{|x|=\delta\}} \Big [ \frac{x_i \, x_j}{r^6} \, D_\nu u - u \, D_\nu \Big ( \frac{x_i \, x_j}{r^6} \Big ) \Big ] \, d\mu_{g_{\text{\rm eucl}}} = \frac{1}{2} \, \pi^2 \, D_i D_j u(0)\]
and 
\[\int_{\{|x|=\delta\}} \Big [ \frac{x_i^2-x_j^2}{r^6} \, D_\nu u - u \, D_\nu \Big ( \frac{x_i^2-x_j^2}{r^6} \Big ) \Big ] \, d\mu_{g_{\text{\rm eucl}}} = \frac{1}{2} \, \pi^2 \, (D_i D_i u(0) - D_j D_j u(0)).\]
Here, $\nu$ denotes the outward-pointing unit normal to the hypersurface $\{|x|=\delta\}$ with respect to the Euclidean metric. 
\end{lemma}

\begin{proof} 
By the divergence theorem, the quantity 
\[I(r) := \int_{\{|x|=r\}} \Big [ \frac{x_i \, x_j}{r^6} \, D_\nu u - u \, D_\nu \Big ( \frac{x_i \, x_j}{r^6} \Big ) \Big ] \, d\mu_{g_{\text{\rm eucl}}}\] 
is independent of $r$. Moreover, if we perform a Taylor expansion of $u$ around the origin, it is easy to see that $I(r) \to \frac{1}{2} \, \pi^2 \, D_i D_j u(0)$ as $r \to 0$. Thus, $I(r) = \frac{1}{2} \, \pi^2 \, D_i D_j u(0)$ for all $r > 0$. This proves the first identity. The second identity follows from an analogous argument.
\end{proof}

\begin{proposition} 
\label{flux.integral}
We have 
\[\int_{\{|x|=\delta\}} (\langle o_{1,\varepsilon},D_\nu \bar{h} \rangle_{g_{\text{\rm eh},\varepsilon}} - \langle \bar{h},D_\nu o_{1,\varepsilon} \rangle_{g_{\text{\rm eh},\varepsilon}}) \, d\mu_{g_{\text{\rm eh},\varepsilon}} = 32\pi^2\omega \, \varepsilon^8 + O(\varepsilon^{12} \, \delta^{-10}),\]
where $\bar{h} = \bar{g}_{\varepsilon,\delta} - g_{\text{\rm eh},\varepsilon}$ and 
\[\omega := \sum_{a \in \mathbb{Z}_{\text{\rm odd}}^4} |a|^{-10} \, (|a|^4 - 6 \, (a_1^2+a_2^2) \, (a_3^2+a_4^2)) \approx 7.70.\]  
\end{proposition}

\begin{proof} 
For $\delta \leq |x| \leq 2\delta$, we have 
\[\bar{h} = \bar{g}_{\varepsilon,\delta} - g_{\text{\rm eh},\varepsilon} = \frac{1}{2} \, \varepsilon^4 \sum_{a \in \mathbb{Z}_{\text{\rm even}}^4 \setminus \{0\}} \tau_a^* T + \frac{1}{2} \, \varepsilon^4 \sum_{a \in \mathbb{Z}_{\text{\rm odd}}^4} \tau_a^* \hat{T} + O(\varepsilon^8 \, |x|^{-8})\] 
and 
\[o_{1,\varepsilon} = \varepsilon^4 \, T + O(\varepsilon^8 \, |x|^{-8}).\] 
Consequently, 
\begin{align*} 
&\int_{\{|x|=\delta\}} (\langle o_{1,\varepsilon},D_\nu \bar{h} \rangle_{g_{\text{\rm eh},\varepsilon}} - \langle \bar{h},D_\nu o_{1,\varepsilon} \rangle_{g_{\text{\rm eh},\varepsilon}}) \, d\mu_{g_{\text{\rm eh},\varepsilon}} \\ 
&= \frac{1}{2} \, \varepsilon^8 \sum_{a \in \mathbb{Z}_{\text{\rm even}}^4 \setminus \{0\}} \int_{\{|x|=\delta\}} (\langle T,D_\nu(\tau_a^* T) \rangle_{g_{\text{\rm eucl}}} - \langle \tau_a^* T,D_\nu T \rangle_{g_{\text{\rm eucl}}}) \, d\mu_{g_{\text{\rm eucl}}} \\ 
&+ \frac{1}{2} \, \varepsilon^8 \sum_{a \in \mathbb{Z}_{\text{\rm odd}}^4} \int_{\{|x|=\delta\}} (\langle T,D_\nu(\tau_a^* \hat{T}) \rangle_{g_{\text{\rm eucl}}} - \langle \tau_a^* \hat{T},D_\nu T \rangle_{g_{\text{\rm eucl}}}) \, d\mu_{g_{\text{\rm eucl}}} \\ 
&+ O(\varepsilon^{12} \, \delta^{-10}).
\end{align*} 
Hence, it remains to evaluate the integrals 
\[\int_{\{|x|=\delta\}} (\langle T,D_\nu(\tau_a^* T) \rangle_{g_{\text{\rm eucl}}} - \langle \tau_a^* T,D_\nu T \rangle_{g_{\text{\rm eucl}}}) \, d\mu_{g_{\text{\rm eucl}}}\] 
and 
\[\int_{\{|x|=\delta\}} (\langle T,D_\nu(\tau_a^* \hat{T}) \rangle_{g_{\text{\rm eucl}}} - \langle \tau_a^* \hat{T},D_\nu T \rangle_{g_{\text{\rm eucl}}}) \, d\mu_{g_{\text{\rm eucl}}},\] 
where $a \in \mathbb{Z}^4 \setminus \{0\}$. To that end, we use Lemma \ref{distributional.laplacian}. Since the components of $\tau_a^* T$ are smooth harmonic functions near the origin, we obtain 
\begin{align*} 
&\int_{\{|x|=\delta\}} (\langle T,D_\nu(\tau_a^* T) \rangle_{g_{\text{\rm eucl}}} - \langle \tau_a^* T,D_\nu T \rangle_{g_{\text{\rm eucl}}}) \, \, d\mu_{g_{\text{\rm eucl}}} \\ 
&= -\frac{1}{2} \, \pi^2 \, \langle (D_1 D_1 + D_2 D_2 - D_3 D_3 - D_4 D_4) \tau_a^* T, \\ 
&\hspace{20mm} dx_1 \otimes dx_1 + dx_2 \otimes dx_2 - dx_3 \otimes dx_3 - dx_4 \otimes dx_4 \rangle \Big |_{x=0} \\ 
&- \pi^2 \, \langle (D_1 D_3 + D_2 D_4) \tau_a^* T,dx_1 \otimes dx_3 + dx_3 \otimes dx_1 + dx_2 \otimes dx_4 + dx_4 \otimes dx_2 \rangle \Big |_{x=0} \\ 
&- \pi^2 \, \langle (D_1 D_4 - D_2 D_3) \tau_a^* T,dx_1 \otimes dx_4 + dx_4 \otimes dx_1 - dx_2 \otimes dx_3 - dx_3 \otimes dx_2 \rangle \Big |_{x=0} \\ 
&= 2\pi^2 \, (D_1 D_1 + D_2 D_2 - D_3 D_3 - D_4 D_4) \\ 
&\hspace{20mm} \frac{(x_1-a_1)^2+(x_2-a_2)^2-(x_3-a_3)^2-(x_4-a_4)^2}{|x-a|^6} \Big |_{x=0} \\ 
&+ 8\pi^2 \, (D_1 D_3 + D_2 D_4) \frac{(x_1-a_1)(x_3-a_3)+(x_2-a_2)(x_4-a_4)}{|x-a|^6} \Big |_{x=0} \\ 
&+ 8\pi^2 \, (D_1 D_4 - D_2 D_3) \frac{(x_1-a_1)(x_4-a_4)-(x_2-a_2)(x_3-a_3)}{|x-a|^6} \Big |_{x=0} \\ 
&= 0. 
\end{align*} 
Moreover, since the components of $\tau_a^* \hat{T}$ are smooth harmonic functions near the origin, we obtain 
\begin{align*} 
&\int_{\{|x|=\delta\}} (\langle T,D_\nu(\tau_a^* \hat{T}) \rangle_{g_{\text{\rm eucl}}} - \langle \tau_a^* \hat{T},D_\nu T \rangle_{g_{\text{\rm eucl}}}) \, d\mu_{g_{\text{\rm eucl}}} \\ 
&= -\frac{1}{2} \, \pi^2 \, \langle (D_1 D_1 + D_2 D_2 - D_3 D_3 - D_4 D_4) \tau_a^* \hat{T}, \\ 
&\hspace{20mm} dx_1 \otimes dx_1 + dx_2 \otimes dx_2 - dx_3 \otimes dx_3 - dx_4 \otimes dx_4 \rangle \Big |_{x=0} \\ 
&- \pi^2 \, \langle (D_1 D_3 + D_2 D_4) \tau_a^* \hat{T},dx_1 \otimes dx_3 + dx_3 \otimes dx_1 + dx_2 \otimes dx_4 + dx_4 \otimes dx_2 \rangle \Big |_{x=0} \\ 
&- \pi^2 \, \langle (D_1 D_4 - D_2 D_3) \tau_a^* \hat{T},dx_1 \otimes dx_4 + dx_4 \otimes dx_1 - dx_2 \otimes dx_3 - dx_3 \otimes dx_2 \rangle \Big |_{x=0} \\ 
&= 2\pi^2 \, (D_1 D_1 + D_2 D_2 - D_3 D_3 - D_4 D_4) \\ 
&\hspace{20mm} \frac{(x_1-a_1)^2+(x_2-a_2)^2-(x_3-a_3)^2-(x_4-a_4)^2}{|x-a|^6} \Big |_{x=0} \\ 
&= 64\pi^2 \, |a|^{-10} \, (|a|^4 - 6 \, (a_1^2+a_2^2) \, (a_3^2+a_4^2)). 
\end{align*} 
Putting these facts together, the assertion follows.
\end{proof}

\begin{proposition}
\label{Z.term}
Let $\bar{h} = \bar{g}_{\varepsilon,\delta} - g_{\text{\rm eh},\varepsilon}$ and $Z = \text{\rm div}_{g_{\text{\rm eh},\varepsilon}} \, \bar{h} - \frac{1}{2} \, \nabla \text{\rm tr}_{g_{\text{\rm eh},\varepsilon}} \, \bar{h}$. Then 
\[\int_{\{|x|=\delta\}} o_{1,\varepsilon}(Z,\nu) \, d\mu_{g_{\text{\rm eh},\varepsilon}} = O(\varepsilon^{12} \, \delta^{-10}).\] 
\end{proposition} 

\begin{proof} 
Recall that
\[\bar{h} = \bar{g}_{\varepsilon,\delta} - g_{\text{\rm eh},\varepsilon} = \frac{1}{2} \, \varepsilon^4 \sum_{a \in \mathbb{Z}_{\text{\rm even}}^4 \setminus \{0\}} \tau_a^* T + \frac{1}{2} \, \varepsilon^4 \sum_{a \in \mathbb{Z}_{\text{\rm odd}}^4} \tau_a^* \hat{T} + O(\varepsilon^8 \, |x|^{-8})\] 
for $\delta \leq |x| \leq 2\delta$. Both $\tau_a^* T$ and $\tau_a^* \hat{T}$ are trace-free and divergence-free with respect to the Euclidean metric. This implies 
\[\nabla \text{\rm tr}_{g_{\text{\rm eucl}}} \bar{h} = O(\varepsilon^8 \, |x|^{-9})\] 
and 
\[\text{\rm div}_{g_{\text{\rm eucl}}} \bar{h} = O(\varepsilon^8 \, |x|^{-9})\] 
for $\delta \leq |x| \leq 2\delta$. Since $g_{\text{\rm eh},\varepsilon} - g_{\text{\rm eucl}} = O(\varepsilon^4 \, |x|^{-4})$, it follows that  
\[\nabla \text{\rm tr}_{g_{\text{\rm eh},\varepsilon}} \bar{h} = O(\varepsilon^8 \, |x|^{-9})\] 
and 
\[\text{\rm div}_{g_{\text{\rm eh},\varepsilon}} \bar{h} = O(\varepsilon^8 \, |x|^{-9})\] 
for $\delta \leq |x| \leq 2\delta$. In particular, $Z = O(\varepsilon^8 \, |x|^{-9})$, which implies $o_{1,\varepsilon}(Z,\nu) = O(\varepsilon^{12} \, \delta^{-13})$. From this, the assertion follows easily.
\end{proof}

\begin{corollary} 
\label{projection.of.Ric.to.o_1}
We have 
\[-2 \int_{[-\frac{1}{2},\frac{1}{2}]^4 \setminus \{0\}} \langle \bar{o}_{1,\varepsilon,\delta},\text{\rm Ric}_{\bar{g}_{\varepsilon,\delta}} \rangle_{\bar{g}_{\varepsilon,\delta}} \, d\text{\rm vol}_{\bar{g}_{\varepsilon,\delta}} = 32\pi^2\omega \, \varepsilon^8 + O(\varepsilon^8 \, \delta^2).\]
\end{corollary}

\begin{proof} 
As above, we write $\bar{h} = \bar{g}_{\varepsilon,\delta} - g_{\text{\rm eh},\varepsilon}$ and $Z = \text{\rm div}_{g_{\text{\rm eh},\varepsilon}} \, \bar{h} - \frac{1}{2} \, \nabla \text{\rm tr}_{g_{\text{\rm eh},\varepsilon}} \, \bar{h}$. Note that $\bar{h}$ is defined in the region $\{|x| \leq \delta\}$ and vanishes identically in the region $\{|x| \leq \frac{1}{2} \, \delta\}$. Moreover, we have $|\bar{h}| \leq O(\varepsilon^4)$, $|\nabla \bar{h}| \leq O(\varepsilon^4 \, \delta^{-1})$, and $|\nabla^2 \bar{h}| \leq O(\varepsilon^4 \, \delta^{-2})$ in the annulus $\{\frac{1}{2} \, \delta \leq |x| \leq \delta\}$. 

Using the well-known formula for the linearization of the Ricci curvature, we obtain 
\[-2 \, \text{\rm Ric}_{\bar{g}_{\varepsilon,\delta}} = \Delta_{L,g_{\text{\rm eh},\varepsilon}} \bar{h} - \mathscr{L}_Z g_{\text{\rm eh},\varepsilon} + O(\varepsilon^8 \, \delta^{-2})\] 
in the region $\{\frac{1}{2} \, \delta \leq |x| \leq \delta\}$. Since $\bar{o}_{1,\varepsilon,\delta} = o_{1,\varepsilon} + O(\varepsilon^4)$ and $\text{\rm Ric}_{\bar{g}_{\varepsilon,\delta}} = O(\varepsilon^4 \, \delta^{-2})$ in the region $\{\frac{1}{2} \, \delta \leq |x| \leq \delta\}$, we deduce that 
\begin{align*} 
&-2 \int_{\{\frac{1}{2} \, \delta \leq |x| \leq \delta\}} \langle \bar{o}_{1,\varepsilon,\delta},\text{\rm Ric}_{\bar{g}_{\varepsilon,\delta}} \rangle_{\bar{g}_{\varepsilon,\delta}} \, d\text{\rm vol}_{\bar{g}_{\varepsilon,\delta}} \\ 
&= \int_{\{\frac{1}{2} \, \delta \leq |x| \leq \delta\}} \langle o_{1,\varepsilon},\Delta_{L,g_{\text{\rm eh},\varepsilon}} \bar{h} - \mathscr{L}_Z g_{\text{\rm eh},\varepsilon} \rangle_{g_{\text{\rm eh},\varepsilon}} \, d\text{\rm vol}_{g_{\text{\rm eh},\varepsilon}} + O(\varepsilon^8 \, \delta^2). 
\end{align*} 
Since $o_{1,\varepsilon}$ satisfies the equation $\Delta_{L,g_{\text{\rm eh},\varepsilon}} o_{1,\varepsilon} = 0$, the divergence theorem gives 
\begin{align*} 
&\int_{\{\frac{1}{2} \, \delta \leq |x| \leq \delta\}} \langle o_{1,\varepsilon},\Delta_{L,g_{\text{\rm eh},\varepsilon}} \bar{h} \rangle_{g_{\text{\rm eh},\varepsilon}} \, d\text{\rm vol}_{g_{\text{\rm eh},\varepsilon}} \\ 
&= \int_{\{|x| \leq \delta\}} (\langle o_{1,\varepsilon},\Delta_{L,g_{\text{\rm eh},\varepsilon}} \bar{h} \rangle_{g_{\text{\rm eh},\varepsilon}} - \langle \bar{h},\Delta_{L,g_{\text{\rm eh},\varepsilon}} o_{1,\varepsilon} \rangle_{g_{\text{\rm eh},\varepsilon}}) \, d\text{\rm vol}_{g_{\text{\rm eh},\varepsilon}} \\ 
&= \int_{\{|x|=\delta\}} (\langle o_{1,\varepsilon},D_\nu \bar{h} \rangle_{g_{\text{\rm eh},\varepsilon}} - \langle \bar{h},D_\nu o_{1,\varepsilon} \rangle_{g_{\text{\rm eh},\varepsilon}}) \, d\mu_{g_{\text{\rm eh},\varepsilon}} \\ 
&= 32\pi^2\omega \, \varepsilon^8 + O(\varepsilon^{12} \, \delta^{-10}), 
\end{align*} 
where in the last step we have used Proposition \ref{flux.integral}. Finally, since $\text{\rm div}_{g_{\text{\rm eh},\varepsilon}} o_{1,\varepsilon} = 0$, the divergence theorem yields 
\begin{align*} 
&\int_{\{\frac{1}{2} \, \delta \leq |x| \leq \delta\}} \langle o_{1,\varepsilon},\mathscr{L}_Z g_{\text{\rm eh},\varepsilon} \rangle_{g_{\text{\rm eh},\varepsilon}} \, d\text{\rm vol}_{g_{\text{\rm eh},\varepsilon}} \\ 
&= \int_{\{|x| \leq \delta\}} (\langle o_{1,\varepsilon},\mathscr{L}_Z g_{\text{\rm eh},\varepsilon} \rangle_{g_{\text{\rm eh},\varepsilon}} + 2 \, \langle \text{\rm div}_{g_{\text{\rm eh},\varepsilon}} o_{1,\varepsilon},Z \rangle_{g_{\text{\rm eh},\varepsilon}}) \, d\text{\rm vol}_{g_{\text{\rm eh},\varepsilon}} \\ 
&= \int_{\{|x|=\delta\}} 2 \, o_{1,\varepsilon}(Z,\nu) \, d\mu_{g_{\text{\rm eh},\varepsilon}} \\ 
&= O(\varepsilon^{12} \, \delta^{-10}) 
\end{align*} 
by Proposition \ref{Z.term}. Putting these facts together, we conclude that 
\[-2 \int_{\{\frac{1}{2} \, \delta \leq |x| \leq \delta\}} \langle \bar{o}_{1,\varepsilon,\delta},\text{\rm Ric}_{\bar{g}_{\varepsilon,\delta}} \rangle_{\bar{g}_{\varepsilon,\delta}} \, d\text{\rm vol}_{\bar{g}_{\varepsilon,\delta}} = 32\pi^2\omega \, \varepsilon^8 + O(\varepsilon^8 \, \delta^2).\] 
Since $\text{\rm Ric}_{\bar{g}_{\varepsilon,\delta}} = 0$ in the region $\{|x| \leq \frac{1}{2} \, \delta\}$ and $|\text{\rm Ric}_{\bar{g}_{\varepsilon,\delta}}| \leq C \, \varepsilon^8 \, |x|^{-10}$ in the region $[-\frac{1}{2},\frac{1}{2}]^4 \setminus \{|x| \leq \delta\}$, the assertion follows.
\end{proof}

\begin{proposition}
\label{projection.of.Ric.to.g}
We have 
\[-2 \displaystyle\int_{[-\frac{1}{2},\frac{1}{2}]^4 \setminus \{0\}} \langle \bar{g}_{\varepsilon,\delta},\text{\rm Ric}_{\bar{g}_{\varepsilon,\delta}} 
\rangle_{\bar{g}_{\varepsilon,\delta}} \, d\text{\rm vol}_{\bar{g}_{\varepsilon,\delta}} = O(\varepsilon^8 \, \delta^{-6}).\] 
\end{proposition}

\begin{proof} 
Let $\tilde{h} = \bar{g}_{\varepsilon,\delta} - g_{\text{\rm eucl}}$. Clearly, $\tilde{h}$ is defined on $[-\frac{1}{2},\frac{1}{2}]^4 \setminus \{|x| \leq \frac{1}{4} \, \delta\}$, and we have $|\tilde{h}| \leq O(\varepsilon^4 \, |x|^{-4})$, $|\nabla \tilde{h}| \leq O(\varepsilon^4 \, |x|^{-5})$, and $|\nabla^2 \tilde{h}| \leq O(\varepsilon^4 \, |x|^{-6})$. Using the standard formula for the linearization of the scalar curvature, we deduce that 
\[\langle \bar{g}_{\varepsilon,\delta},\text{\rm Ric}_{\bar{g}_{\varepsilon,\delta}} \rangle_{\bar{g}_{\varepsilon,\delta}} = \text{\rm div} \, \text{\rm div} \, \tilde{h} - \Delta \text{\rm tr} \, \tilde{h} + O(\varepsilon^8 \, |x|^{-10})\] 
in the region $[-\frac{1}{2},\frac{1}{2}]^4 \setminus \{|x| \leq \frac{1}{2} \, \delta\}$. Here, the divergence and the Laplacian on the right hand side are taken with respect to the Euclidean metric. Using the divergence theorem, we obtain 
\begin{align*} 
&\int_{[-\frac{1}{2},\frac{1}{2}]^4 \setminus \{|x| \leq \frac{1}{2} \, \delta\}} \langle \bar{g}_{\varepsilon,\delta},\text{\rm Ric}_{\bar{g}_{\varepsilon,\delta}} \rangle_{\bar{g}_{\varepsilon,\delta}} \, d\text{\rm vol}_{\bar{g}_{\varepsilon,\delta}} \\ 
&= \int_{[-\frac{1}{2},\frac{1}{2}]^4 \setminus \{|x| \leq \frac{1}{2} \, \delta\}} (\text{\rm div} \, \text{\rm div} \, \tilde{h} - \Delta \text{\rm tr} \, \tilde{h}) \, d\text{\rm vol}_{g_{\text{\rm eucl}}} + O(\varepsilon^8 \, \delta^{-6}) \\ 
&= -\int_{\{|x|=\frac{1}{2} \, \delta\}} (\langle \text{\rm div} \, \tilde{h},\nu \rangle_{g_{\text{\rm eucl}}} - \langle \nabla \text{\rm tr} \, \tilde{h},\nu \rangle_{g_{\text{\rm eucl}}}) \, d\mu_{g_{\text{\rm eucl}}} + O(\varepsilon^8 \, \delta^{-6}). 
\end{align*}
As above, $\nu$ denotes the outward-pointing unit normal to the hypersurface $\{|x|=\frac{1}{2} \, \delta\}$ with respect to the Euclidean metric. For $\frac{1}{4} \, \delta \leq |x| \leq \frac{1}{2} \, \delta$, we have 
\[\tilde{h} = g_{\text{\rm eh},\varepsilon} - g_{\text{\rm eucl}} = \frac{1}{2} \, \varepsilon^4 \, T + O(\varepsilon^8 \, \delta^{-8}).\] 
Since $T$ is trace-free and divergence-free with respect to the Euclidean metric, it follows that 
\[\nabla \text{\rm tr} \, \tilde{h} = O(\varepsilon^8 \, \delta^{-9})\] 
and 
\[\text{\rm div} \, \tilde{h} = O(\varepsilon^8 \, \delta^{-9}).\] 
Thus, we conclude that 
\[\int_{\{|x|=\frac{1}{2} \, \delta\}} (\langle \text{\rm div} \, \tilde{h},\nu \rangle_{g_{\text{\rm eucl}}} - \langle \nabla \text{\rm tr} \, \tilde{h},\nu \rangle_{g_{\text{\rm eucl}}}) \, d\mu_{g_{\text{\rm eucl}}} = O(\varepsilon^8 \, \delta^{-6}),\] 
hence 
\[\int_{[-\frac{1}{2},\frac{1}{2}]^4 \setminus \{|x| \leq \frac{1}{2} \, \delta\}} \langle \bar{g}_{\varepsilon,\delta},\text{\rm Ric}_{\bar{g}_{\varepsilon,\delta}} \rangle_{\bar{g}_{\varepsilon,\delta}} \, d\text{\rm vol}_{\bar{g}_{\varepsilon,\delta}} = O(\varepsilon^8 \, \delta^{-6}).\] 
Since $\text{\rm Ric}_{\bar{g}_{\varepsilon,\delta}} = 0$ in the region $\{|x| \leq \frac{1}{2} \, \delta\}$, the assertion follows.
\end{proof}

Finally, we compute the projection of $\frac{\partial}{\partial \varepsilon} \bar{g}_{\varepsilon,\delta}$ to the approximate kernel. 

\begin{proposition} 
\label{projection.of.time.derivative}
We have 
\[\int_{[-\frac{1}{2},\frac{1}{2}]^4 \setminus \{0\}} \langle \bar{o}_{1,\varepsilon,\delta},\frac{\partial}{\partial \varepsilon} \bar{g}_{\varepsilon,\delta} \rangle_{\bar{g}_{\varepsilon,\delta}} \, d\text{\rm vol}_{\bar{g}_{\varepsilon,\delta}} = 4\pi^2 \varepsilon^3 + O(\varepsilon^7 \delta^{-4})\] 
and 
\[\int_{[-\frac{1}{2},\frac{1}{2}]^4 \setminus \{0\}} \langle \bar{g}_{\varepsilon,\delta},\frac{\partial}{\partial \varepsilon} \bar{g}_{\varepsilon,\delta} \rangle_{\bar{g}_{\varepsilon,\delta}} \, d\text{\rm vol}_{\bar{g}_{\varepsilon,\delta}} = O(\varepsilon^3 \log \frac{1}{\delta}).\] 
\end{proposition} 

\begin{proof} 
For $|x| \leq \frac{1}{2} \, \delta$, we have $\bar{g}_{\varepsilon,\delta} = g_{\text{\rm eh},\varepsilon}$, $\frac{\partial}{\partial \varepsilon} \bar{g}_{\varepsilon,\delta} = \frac{2}{\varepsilon} \, o_{1,\varepsilon}$, and $\bar{o}_{1,\varepsilon,\delta} = o_{1,\varepsilon}$. From this, the assertion follows easily.
\end{proof}

\section{Liouville-type theorems for the parabolic Lichnerowicz equation}

\label{liouville}

In this section, we establish Liouville-type theorems for the parabolic Lichnerowicz equation on various model spaces. We first recall some basic properties of the Lichnerowicz Laplacian on the Eguchi-Hanson manifold which were established by Biquard \cite{Biquard} and Biquard and Rollin \cite{Biquard-Rollin} (see also \cite{Page1}). 

\begin{proposition}[O.~Biquard \cite{Biquard}; O.~Biquard, Y.~Rollin \cite{Biquard-Rollin}]
\label{kernel}
Let $(M_{\text{\rm eh}},g_{\text{\rm eh}})$ denote the Eguchi-Hanson manifold with parameter $\varepsilon=1$, and let $o_i := o_{i,1}$. Then the following statements hold: 
\begin{itemize}
\item[(i)] Let $h$ be a symmetric $(0,2)$-tensor on the Eguchi-Hanson manifold satisfying $|h| \leq (1+r)^{-\sigma}$ for some $\sigma > 0$ and $\Delta_{L,g_{\text{\rm eh}}} h = 0$. Then $h \in \text{\rm span}\{o_1,o_2,o_3\}$. 
\item[(ii)] Let $h$ be a symmetric $(0,2)$-tensor on the Eguchi-Hanson manifold satisfying $|h| \leq (1+r)^{-\sigma-1}$ and $|\nabla h| \leq (1+r)^{-\sigma-2}$ for some $\sigma > 0$. Then 
\[\int_{M_{\text{\rm eh}}} \langle \Delta_{L,g_{\text{\rm eh}}} h,h \rangle \, d\text{\rm vol}_{g_{\text{\rm eh}}} \leq 0.\] 
Moreover, equality holds if and only if $h \in \text{\rm span} \{o_1,o_2,o_3\}$.
\end{itemize}
\end{proposition}

\begin{proof}
The first statement is contained in Proposition 1.1 in \cite{Biquard}. The second statement follows from a Bochner formula from \cite{Biquard-Rollin}. To explain this, let $h$ be a symmetric $(0,2)$-tensor on the Eguchi-Hanson manifold satisfying $|h| \leq C \, (1+r)^{-\sigma-1}$ and $|\nabla h| \leq C \, (1+r)^{-\sigma-2}$. Clearly, 
\[\int_{M_{\text{\rm eh}}} \langle \Delta_{L,g_{\text{\rm eh}}} h,h \rangle \, d\text{\rm vol}_{g_{\text{\rm eh}}} = \int_{M_{\text{\rm eh}}} \langle \Delta_{L,g_{\text{\rm eh}}} k,k \rangle \, d\text{\rm vol}_{g_{\text{\rm eh}}} + \frac{1}{4} \int_{M_{\text{\rm eh}}} \langle \Delta_{g_{\text{\rm eh}}} \text{\rm tr} \, h,\text{\rm tr} \, h \rangle \, d\text{\rm vol}_{g_{\text{\rm eh}}},\]
where $k$ denotes the trace-free part of $h$. We may view $k$ as a section of the vector bundle $\Lambda_-^2 \otimes \Lambda_+^2$. The Bochner formula (4.6) in \cite{Biquard-Rollin} implies that $-\frac{1}{2} \, \Delta_{L,g_{\text{\rm eh}}} k = d_- d_-^* k$, where $d_-: \Gamma(\Lambda^1 \otimes \Lambda_+^2) \to \Gamma(\Lambda_-^2 \otimes \Lambda_+^2)$ is the exterior derivative (cf. \cite{Biquard}, \cite{Biquard-Rollin}). This gives 
\[\int_{M_{\text{\rm eh}}} \langle \Delta_{L,g_{\text{\rm eh}}} k,k \rangle \, d\text{\rm vol}_{g_{\text{\rm eh}}} = -2 \int_{M_{\text{\rm eh}}} |d_-^* k|^2 \, d\text{\rm vol}_{g_{\text{\rm eh}}},\] 
hence 
\[\int_{M_{\text{\rm eh}}} \langle \Delta_{L,g_{\text{\rm eh}}} h,h \rangle \, d\text{\rm vol}_{g_{\text{\rm eh}}} = -2 \int_{M_{\text{\rm eh}}} |d_-^* k|^2 \, d\text{\rm vol}_{g_{\text{\rm eh}}} - \frac{1}{4} \int_{M_{\text{\rm eh}}} |\nabla \text{\rm tr} \, h|^2 \, d\text{\rm vol}_{g_{\text{\rm eh}}} \leq 0.\] 
Moreover, if equality holds, then $d_-^* k = 0$ and $\text{\rm tr} \, h = 0$. Results in \cite{Biquard} now imply that $h \in \text{\rm span}\{o_1,o_2,o_3\}$.
\end{proof}

We next establish a Liouville-type theorem for the linear heat equation on the Eguchi-Hanson manifold.

\begin{proposition}
\label{liouville.thm.1}
Let $(M_{\text{\rm eh}},g_{\text{\rm eh}})$ denote the Eguchi-Hanson manifold with parameter $\varepsilon=1$, and let $o_i := o_{i,1}$. Let $h$ be a solution of the heat equation $\frac{\partial}{\partial t} h = \Delta_{L,g_{\text{\rm eh}}} h$ on $M_{\text{\rm eh}} \times (-\infty,0]$ with the property that $|h| \leq (1+r)^{-\sigma}$ for some $\sigma>0$. If $\int_{M_{\text{\rm eh}}} \langle h(t),o_i \rangle = 0$ for all $i \in \{1,2,3\}$ and all $t \in (-\infty,0]$, then $h$ vanishes identically.
\end{proposition}

\begin{proof} 
Let us consider the tensor field $k(t) := \Delta_{L,g_{\text{\rm eh}}} h(t)$. It follows from standard interior estimates for parabolic equations that $|k| \leq C \, (1+r)^{-\sigma-2}$ and $|\nabla k| \leq C \, (1+r)^{-\sigma-3}$. In particular, $\int_{M_{\text{\rm eh}}} |k(t)|^2 \leq C$ for each $t \in (-\infty,0]$. Moreover, $k$ satisfies the equation $\frac{\partial}{\partial t} k = \Delta_{L,g_{\text{\rm eh}}} k$, and we have $\int_{M_{\text{\rm eh}}} \langle k(t),o_i \rangle = 0$ for all $i \in \{1,2,3\}$ and all $t \in (-\infty,0]$. Using Proposition \ref{kernel}, we obtain 
\[\frac{1}{2} \, \frac{d}{dt} \bigg ( \int_{M_{\text{\rm eh}}} |k(t)|^2 \, d\text{\rm vol}_{g_{\text{\rm eh}}} \bigg ) = \int_{M_{\text{\rm eh}}} \langle \Delta_{L,g_{\text{\rm eh}}} k(t),k(t) \rangle \, d\text{\rm vol}_{g_{\text{\rm eh}}} \leq 0.\] 
Consequently, the function $t \mapsto \int_{M_{\text{\rm eh}}} |k(t)|^2$ is monotone decreasing. In particular, the limit $A := \lim_{t \to -\infty} \int_{M_{\text{\rm eh}}} |k(t)|^2$ exists.

We next pick an arbitrary sequence of times $t_j \to -\infty$, and define $\tilde{k}^{(j)}(t) := k(t_j+t)$. After passing to a subsequence if necessary, we may assume that the sequence $\tilde{k}^{(j)}$ converges in $C_{loc}^\infty$ to some tensor field $\hat{k}$ which is defined on $M_{\text{\rm eh}} \times \mathbb{R}$ and satisfies the equation $\frac{\partial}{\partial t} \hat{k} = \Delta_{L,g_{\text{\rm eh}}} \hat{k}$. Moreover, $|\hat{k}| \leq C \, (1+r)^{-\sigma-2}$ and $|\nabla \hat{k}| \leq C \, (1+r)^{-\sigma-3}$. Using the dominated convergence theorem, we obtain 
\[\int_{M_{\text{\rm eh}}} |\hat{k}(t)|^2 \, d\text{\rm vol}_{g_{\text{\rm eh}}} = \lim_{j \to \infty} \int_{M_{\text{\rm eh}}} |k(t_j+t)|^2 \, d\text{\rm vol}_{g_{\text{\rm eh}}} = A\] 
and 
\[\int_{M_{\text{\rm eh}}} \langle \hat{k}(t),o_i \rangle \, d\text{\rm vol}_{g_{\text{\rm eh}}} = \lim_{j \to \infty} \int_{M_{\text{\rm eh}}} \langle k(t_j+t),o_i \rangle \, d\text{\rm vol}_{g_{\text{\rm eh}}} = 0\] 
for all $i \in \{1,2,3\}$ and all $t \in \mathbb{R}$. Differentiating the first identity with respect to $t$ gives 
\[0 = \frac{1}{2} \, \frac{d}{dt} \bigg ( \int_{M_{\text{\rm eh}}} |\hat{k}(t)|^2 \, d\text{\rm vol}_{g_{\text{\rm eh}}} \bigg ) = \int_{M_{\text{\rm eh}}} \langle \Delta_{L,g_{\text{\rm eh}}} \hat{k}(t),\hat{k}(t) \rangle \, d\text{\rm vol}_{g_{\text{\rm eh}}}\] 
for all $t \in \mathbb{R}$. Hence, Proposition \ref{kernel} implies that $\hat{k}(t) = 0$ for all $t \in \mathbb{R}$. Consequently, $A=0$. In other words, $\lim_{t \to -\infty} \int_{M_{\text{\rm eh}}} |k(t)|^2 = 0$. Since the function $t \mapsto \int_{M_{\text{\rm eh}}} |k(t)|^2$ is monotone decreasing, it follows that $k(t) = 0$ for all $t \in (-\infty,0]$. Thus, we conclude that $\Delta_{L,g_{\text{\rm eh}}} h(t) = 0$ for all $t \in (-\infty,0]$. Since $\int_{M_{\text{\rm eh}}} \langle h(t),o_i \rangle = 0$ for all $i \in \{1,2,3\}$ and all $t \in (-\infty,0]$, Proposition \ref{kernel} implies that $h(t) = 0$ for all $t \in (-\infty,0]$.
\end{proof}

\begin{proposition}
\label{liouville.thm.2}
Let $h$ be a solution of the heat equation $\frac{\partial}{\partial t} h = \Delta_{g_{\text{\rm eucl}}} h$ on $(\mathbb{R}^4 \setminus \{0\}) \times (-\infty,0]$ with the property that $|h| \leq r^{-\sigma}$ for some $\sigma \in (0,2)$. Then $h$ vanishes identically.
\end{proposition}

\begin{proof} 
Let us fix an arbitrary point $(x_0,t_0)$ in spacetime. Since $\sigma < 2$, the equation $\frac{\partial}{\partial t} h = \Delta_{g_{\text{\rm eucl}}} h$ is satisfied in the sense of distributions. Consequently, 
\[h(x_0,t_0) = \int_{\mathbb{R}^4 \setminus \{0\}} \frac{1}{(4\pi t)^2} \, e^{-\frac{|x-x_0|^2}{4t}} \, h(x,t_0-t)\] 
for all $t > 0$. Since $|h(x,t_0-t)| \leq |x|^{-\sigma}$, it follows that 
\begin{align*} 
|h(x_0,t_0)| 
&\leq \int_{\mathbb{R}^4 \setminus \{0\}} \frac{1}{(4\pi t)^2} \, e^{-\frac{|x-x_0|^2}{4t}} \, |x|^{-\sigma} \\ 
&\leq C \, t^{-2} \int_{\{|x|^2 \leq t\}} |x|^{-\sigma} + C \int_{\{|x|^2 \geq t\}} |x|^{-4-\sigma} \\ 
&\leq C \, t^{-\frac{\sigma}{2}} 
\end{align*}
for all $t>0$. Sending $t \to \infty$ gives $h(x_0,t_0) = 0$.
\end{proof}

\begin{proposition}
\label{liouville.thm.3}
Let $h$ be a symmetric $(0,2)$-tensor defined on $(\mathbb{R}^4 \setminus \mathbb{Z}^4) \times (-\infty,0]$ which evolves by the heat equation $\frac{\partial}{\partial t} h = \Delta_{g_{\text{\rm eucl}}} h$ and is invariant under the maps 
\begin{align*} 
&(x_1,x_2,x_3,x_4) \mapsto (1-x_1,x_2,x_3,x_4), \\ 
&(x_1,x_2,x_3,x_4) \mapsto (x_1,1-x_2,x_3,x_4), \\ 
&(x_1,x_2,x_3,x_4) \mapsto (x_1,x_2,1-x_3,x_4), \\ 
&(x_1,x_2,x_3,x_4) \mapsto (x_1,x_2,x_3,1-x_4), \\ 
&(x_1,x_2,x_3,x_4) \mapsto (x_1+2,x_2,x_3,x_4), \\ 
&(x_1,x_2,x_3,x_4) \mapsto (x_1,x_2+2,x_3,x_4), \\ 
&(x_1,x_2,x_3,x_4) \mapsto (x_1,x_2,x_3+2,x_4), \\ 
&(x_1,x_2,x_3,x_4) \mapsto (x_1,x_2,x_3,x_4+2).
\end{align*}
Moreover, we assume that $|h| \leq r^{-\sigma}$ for some $\sigma \in (0,2)$. If $\int_{[-\frac{1}{2},\frac{1}{2}]^4 \setminus \{0\}} h_{ii} = 0$ for each $i \in \{1,2,3,4\}$, then $h$ vanishes identically.
\end{proposition}

\begin{proof} 
Let us fix an arbitrary point $(x_0,t_0)$ in spacetime. Moreover, let  
\[\Gamma_+(x,t) = \frac{1}{(4\pi t)^2} \, \bigg ( \sum_{a \in \mathbb{Z}_{\text{\rm even}}^4} e^{-\frac{|x-x_0-a|^2}{4t}} + \sum_{a \in \mathbb{Z}_{\text{\rm odd}}^4} e^{-\frac{|x-x_0-a|^2}{4t}} \bigg )\] 
and 
\[\Gamma_-(x,t) = \frac{1}{(4\pi t)^2} \, \bigg ( \sum_{a \in \mathbb{Z}_{\text{\rm even}}^4} e^{-\frac{|x-x_0-a|^2}{4t}} - \sum_{a \in \mathbb{Z}_{\text{\rm odd}}^4} e^{-\frac{|x-x_0-a|^2}{4t}} \bigg ).\] 
We first consider an index $i \in \{1,2,3,4\}$. Since $\sigma<2$, $h_{ii}$ satisfies the heat equation in the sense of distributions. Hence, we may write 
\[h_{ii}(x_0,t_0) = \int_{[-\frac{1}{2},\frac{1}{2}]^4 \setminus \{0\}} \Gamma_+(x,t) \, h_{ii}(x,t_0-t) = \int_{[-\frac{1}{2},\frac{1}{2}]^4 \setminus \{0\}} (\Gamma_+(x,t) - 1) \, h_{ii}(x,t_0-t)\] 
for all $t>0$. Since $|h_{ii}(x,t_0-t)| \leq |x|^{-\sigma}$, it follows that  
\[|h_{ii}(x_0,t_0)| \leq \sup_{x \in [-\frac{1}{2},\frac{1}{2}]^4} |\Gamma_+(x,t)-1| \int_{[-\frac{1}{2},\frac{1}{2}]^4 \setminus \{0\}} |x|^{-\sigma}\] 
for all $t>0$. It is well known that the heat kernel on a torus converges to a constant at an exponential rate. This means that $\sup_{x \in [-\frac{1}{2},\frac{1}{2}]^4} |\Gamma_+(x,t)-1|$ converges to $0$ as $t \to \infty$. Consequently, $h_{ii}(x_0,t_0) = 0$.

We next consider a pair of indices $i \neq j$. Since $\sigma<2$, $h_{ij}$ satisfies the heat equation in the sense of distributions. This gives 
\[h_{ij}(x_0,t_0) = \int_{[-\frac{1}{2},\frac{1}{2}]^4 \setminus \{0\}} \Gamma_-(x,t) \, h_{ij}(x,t_0-t)\] 
for all $t>0$. Since $|h_{ij}(x,t_0-t)| \leq |x|^{-\sigma}$, we obtain 
\[|h_{ij}(x_0,t_0)| \leq \sup_{x \in [-\frac{1}{2},\frac{1}{2}]^4} |\Gamma_-(x,t)| \int_{[-\frac{1}{2},\frac{1}{2}]^4 \setminus \{0\}} |x|^{-\sigma}\] 
for all $t>0$. Since $\sup_{x \in [-\frac{1}{2},\frac{1}{2}]^4} |\Gamma_-(x,t)|$ converges to $0$ as $t \to \infty$, we conclude that $h_{ij}(x_0,t_0) = 0$. This completes the proof. 
\end{proof}

\section{Uniform estimates for the linearized equation}

Throughout this section, we will fix a real number $\alpha \in (0,1)$. Let $M$ denote the quotient manifold introduced in Section \ref{attaching}. We will consider a one-parameter family of metrics $\bar{g}_{\varepsilon(t),\delta(t)}$, $t \in (-\infty,-\Lambda]$, where $\delta(t) = (-t)^{-\frac{1}{400}}$. Moreover, the parameter $\varepsilon(t)$ is assumed to satisfy the following conditions: 

\begin{assumption}
\label{assumption.on.epsilon}
The function $\varepsilon(t)$ satisfies $(-64\omega t)^{-\frac{1}{4}} \leq \varepsilon(t) \leq (-16\omega t)^{-\frac{1}{4}}$, 
\\ $\big | \frac{d\varepsilon}{dt}(t) \big | \leq (-t)^{-\frac{5}{4}}$, and $|t-t'|^{-\alpha} \, \big | \frac{d\varepsilon}{dt}(t)-\frac{d\varepsilon}{dt}(t') \big | \leq (-t)^{-\frac{5}{4}+\frac{\alpha}{2}}$ for $0 < |t-t'| \leq (-t)^{-\frac{1}{2}}$. 
\end{assumption}

We first prove a weighted sup-estimate for solutions of the parabolic Lichnerowicz equation. As in \cite{Mazzeo-Pacard1}, we use a blow-up argument together with the Liouville-type theorems established in Section \ref{liouville}.

\begin{proposition} 
\label{weighted.C0.estimate}
Given real numbers $\gamma>0$ and $\sigma \in (0,2)$, we can find real numbers $\Lambda>0$ and $C>0$ with the following property. Suppose that $\varepsilon(t)$ is a function which is defined on the interval $(-\infty,-\Lambda]$ and satisfies Assumption \ref{assumption.on.epsilon}, and let $\delta(t) = (-t)^{-\frac{1}{400}}$. Let $h$ be a solution of the inhomogeneous heat equation $\frac{\partial}{\partial t} h(t) = \Delta_{L,\bar{g}_{\varepsilon(t),\delta(t)}} h(t) + \psi(t)$ which is defined on $M \times (-\infty,-\Lambda]$ and satisfies $\sup_{M \times (-\infty,-\Lambda]} (-t)^\gamma \, (\varepsilon(t)+r)^\sigma \, |h(t)|_{\bar{g}_{\varepsilon(t),\delta(t)}} < \infty$. Moreover, we assume that $h$ is invariant under the group $\mathscr{G}$ and satisfies the orthogonality conditions 
\[\int_{[-\frac{1}{2},\frac{1}{2}]^4 \setminus \{0\}} \langle h(t),\bar{o}_{1,\varepsilon(t),\delta(t)} \rangle_{\bar{g}_{\varepsilon(t),\delta(t)}} \, d\text{\rm vol}_{\bar{g}_{\varepsilon(t),\delta(t)}} = 0\] 
and 
\[\int_{[-\frac{1}{2},\frac{1}{2}]^4 \setminus \{0\}} \langle h(t),\bar{g}_{\varepsilon(t),\delta(t)} \rangle_{\bar{g}_{\varepsilon(t),\delta(t)}} \, d\text{\rm vol}_{\bar{g}_{\varepsilon(t),\delta(t)}} = 0\] 
for all $t \in (-\infty,-\Lambda]$. Then 
\begin{multline*} 
\sup_{M \times (-\infty,-\Lambda]} (-t)^\gamma \, ((-t)^{-\frac{1}{4}}+r)^\sigma \, |h(t)|_{\bar{g}_{\varepsilon(t),\delta(t)}} \leq 
\\ 
\leq C \, \sup_{M \times (-\infty,-\Lambda]} (-t)^\gamma \, ((-t)^{-\frac{1}{4}}+r)^{\sigma+2} \, |\psi(t)|_{\bar{g}_{\varepsilon(t),\delta(t)}}. 
\end{multline*} 
\end{proposition} 

\begin{proof} 
We argue by contradiction. If the assertion is false, we can find a sequence of functions $\varepsilon^{(j)}(t)$ and sequences of tensor fields $h^{(j)}$ and $\psi^{(j)}$ with the following properties: 
\begin{itemize} 
\item The functions $\varepsilon^{(j)}(t)$ are defined on the interval $(-\infty,-j]$ and satisfy Assumption \ref{assumption.on.epsilon}.
\item The tensor fields $h^{(j)}$ and $\psi^{(j)}$ are defined on $M \times (-\infty,-j]$ and satisfy the equation 
\[\frac{\partial}{\partial t} h^{(j)}(t) = \Delta_{L,g^{(j)}(t)} h^{(j)}(t) + \psi^{(j)}(t),\] 
where $g^{(j)}(t) := \bar{g}_{\varepsilon^{(j)}(t),\delta(t)}$.
\item We have 
\[\sup_{M \times (-\infty,-j]} (-t)^\gamma \, ((-t)^{-\frac{1}{4}}+r)^\sigma \, |h^{(j)}(t)|_{g^{(j)}(t)} = 1\] 
and 
\[\sup_{M \times (-\infty,-j]} (-t)^\gamma \, ((-t)^{-\frac{1}{4}}+r)^{\sigma+2} \, |\psi^{(j)}(t)|_{g^{(j)}(t)} \to 0\] 
as $j \to \infty$.
\item The tensor $h^{(j)}$ satisfies the orthogonality conditions 
\[\int_{[-\frac{1}{2},\frac{1}{2}]^4 \setminus \{0\}} \langle h^{(j)}(t),\bar{o}_{1,\varepsilon^{(j)}(t),\delta(t)} \rangle_{g^{(j)}(t)} \, d\text{\rm vol}_{g^{(j)}(t)} = 0\] 
and 
\[\int_{[-\frac{1}{2},\frac{1}{2}]^4 \setminus \{0\}} \langle h^{(j)}(t),g^{(j)}(t) \rangle_{g^{(j)}(t)} \, d\text{\rm vol}_{g^{(j)}(t)} = 0\] 
for all $t \in (-\infty,-j]$.
\end{itemize}
For each $j$, we can pick a point $(x_j,t_j) \in M \times (-\infty,-j]$ such that $(-t_j)^\gamma \, ((-t_j)^{-\frac{1}{4}}+r)^\sigma \, |h^{(j)}(t_j)|_{g^{(j)}(t)} \geq \frac{1}{2}$ at the point $x_j$. After passing to a subsequence, we are in one of the following three cases: 

\textit{Case 1:} Suppose $\lim_{j \to \infty} (-t_j)^{\frac{1}{4}} \, |x_j| < \infty$. Let us consider the rescaled metrics 
\[\tilde{g}^{(j)}(t) := \varepsilon^{(j)}(t_j)^{-2} \, \rho_j^* g^{(j)}(t_j+\varepsilon^{(j)}(t_j)^2 \, t),\] 
where $t \in (-\infty,0]$ and $\rho_j: M_{\text{\rm eh}} \to M_{\text{\rm eh}}$ denotes a dilation in space by the factor $\varepsilon^{(j)}(t_j)$. Moreover, we define 
\[\tilde{h}^{(j)}(t) := (-t_j)^\gamma \, \varepsilon^{(j)}(t_j)^{\sigma-2} \, \rho_j^* h^{(j)}(t_j+\varepsilon^{(j)}(t_j)^2 \, t)\] 
and 
\[\tilde{\psi}^{(j)}(t) := (-t_j)^\gamma \, \varepsilon^{(j)}(t_j)^\sigma \, \rho_j^* \psi^{(j)}(t_j+\varepsilon^{(j)}(t_j)^2 \, t)\] 
for $t \in (-\infty,0]$. Clearly, $\frac{\partial}{\partial t} \tilde{h}^{(j)}(t) = \Delta_{L,\tilde{g}^{(j)}(t)} \tilde{h}^{(j)}(t) + \tilde{\psi}^{(j)}(t)$. After passing to a subsequence, the tensor fields $\tilde{h}^{(j)}$ converge in $C_{loc}^0$ to some tensor field $\hat{h} \neq 0$. The limiting tensor field $\hat{h}$ is defined on $M_{\text{\rm eh}} \times (-\infty,0]$, where $M_{\text{\rm eh}}$ denotes the Eguchi-Hanson manifold with parameter $\varepsilon=1$. Moreover, $\hat{h}$ satisfies the heat equation $\frac{\partial}{\partial t} \hat{h} = \Delta_{L,g_{\text{\rm eh}}} \hat{h}$, 
and we have 
$\sup_{M_{\text{\rm eh}} \times (-\infty,0]} (1+r)^\sigma \, |\hat{h}(t)|_{g_{\text{\rm eh}}} \leq 1$. 

In the next step, we show that $\hat{h}(t)$ is orthogonal to $o_1$ for each $t \in (-\infty,0]$. Indeed, if we fix a number $t \in (-\infty,0]$, then we have 
\[\int_{[-\frac{1}{2},\frac{1}{2}]^4 \setminus \{0\}} \langle h^{(j)}(t_j'),\bar{o}_{1,\varepsilon^{(j)}(t_j'),\delta(t_j')} \rangle_{\bar{g}_{\varepsilon^{(j)}(t_j'),\delta(t_j')}} \, d\text{\rm vol}_{\bar{g}_{\varepsilon^{(j)}(t_j'),\delta(t_j')}} = 0,\] 
where $t_j' := t_j+\varepsilon^{(j)}(t_j)^2 \, t$. Passing to the limit as $j \to \infty$, we obtain 
\[\int_{M_{\text{\rm eh}}} \langle \hat{h}(t),o_1 \rangle_{g_{\text{\rm eh}}} \, d\text{\rm vol}_{g_{\text{\rm eh}}} = 0\] 
by the dominated convergence theorem. 

Finally, since $\hat{h}(t)$ is invariant under the map $(x_1,x_2,x_3,x_4) \mapsto (x_3,x_4,x_1,x_2)$, we have 
\[\int_{M_{\text{\rm eh}}} \langle \hat{h}(t),o_2 \rangle_{g_{\text{\rm eh}}} \, d\text{\rm vol}_{g_{\text{\rm eh}}} = \int_{M_{\text{\rm eh}}} \langle \hat{h}(t),o_3 \rangle_{g_{\text{\rm eh}}} \, d\text{\rm vol}_{g_{\text{\rm eh}}} = 0\] 
for each $t \in (-\infty,0]$. This contradicts Proposition \ref{liouville.thm.1}. 

\textit{Case 2:} Suppose now that $\lim_{j \to \infty} (-t_j)^{\frac{1}{4}} \, |x_j| = \infty$ and $\lim_{j \to \infty} |x_j| = 0$. We consider the rescaled metrics 
\[\tilde{g}^{(j)}(t) := |x_j|^{-2} \, \rho_j^* g^{(j)}(t_j+|x_j|^2 \, t),\] 
where $t \in (-\infty,0]$ and $\rho_j$ denotes a dilation in space by the factor $|x_j|$. Moreover, we define 
\[\tilde{h}^{(j)}(t) := (-t_j)^\gamma \, |x_j|^{\sigma-2} \, \rho_j^* h(t_j+|x_j|^2 \, t)\] 
and 
\[\tilde{\psi}^{(j)}(t) := (-t_j)^\gamma \, |x_j|^\sigma \, \rho_j^* \psi^{(j)}(t_j+|x_j|^2 \, t)\] 
for $t \in (-\infty,0]$. Clearly, $\frac{\partial}{\partial t} \tilde{h}^{(j)}(t) = \Delta_{L,\tilde{g}^{(j)}(t)} \tilde{h}^{(j)}(t) + \tilde{\psi}^{(j)}(t)$. After passing to a subsequence, the tensors $\tilde{h}^{(j)}$ converge in $C_{loc}^0$ to a tensor field $\hat{h} \neq 0$. The tensor field $\hat{h}$ is defined on $(\mathbb{R}^4 \setminus \{0\}) / \mathbb{Z}_2 \times (-\infty,0]$ and satisfies the equation $\frac{\partial}{\partial t} \hat{h} = \Delta_{g_{\text{\rm eucl}}} \hat{h}$. Moreover, we have $\sup_{(\mathbb{R}^4 \setminus \{0\}) / \mathbb{Z}_2 \times (-\infty,0]} r^\sigma \, |\hat{h}(t)|_{g_{\text{\rm eucl}}} \leq 1$. After lifting $\hat{h}$ to a solution of the heat equation on $(\mathbb{R}^4 \setminus \{0\}) \times (-\infty,0]$, we obtain a contradiction with Proposition \ref{liouville.thm.2}. 

\textit{Case 3:} Suppose finally that $\lim_{j \to \infty} |x_j| > 0$. In this case, we define \[\tilde{h}^{(j)}(t) := (-t_j)^\gamma \, h(t_j+t)\] 
for $t \in (-\infty,0]$. After passing to a subsequence, the tensor fields $\tilde{h}^{(j)}$ converge in $C_{loc}^0$ to some tensor field $\hat{h} \neq 0$. After passing to a suitable covering, we may view $\hat{h}$ as a tensor field on $(\mathbb{R}^4 \setminus \mathbb{Z}^4) \times (-\infty,0]$ which is invariant under the group $\mathscr{G}$. Moreover, $\hat{h}$ satisfies the heat equation $\frac{\partial}{\partial t} \hat{h} = \Delta_{g_{\text{\rm eucl}}} \hat{h}$, 
and we have \\ 
$\sup_{([-\frac{1}{2},\frac{1}{2}]^4 \setminus \{0\}) \times (-\infty,0]} r^\sigma \, |\hat{h}(t)|_{g_{\text{\rm eucl}}} \leq 1$. 

In the next step, we show that $\hat{h}(t)$ is orthogonal to $g_{\text{\rm eucl}}$ for each $t \in (-\infty,0]$. Indeed, if we fix a number $t \in (-\infty,0]$, then we have 
\[\int_{[-\frac{1}{2},\frac{1}{2}]^4 \setminus \{0\}} \langle h^{(j)}(t_j'),\bar{g}_{\varepsilon^{(j)}(t_j'),\delta(t_j')} \rangle_{\bar{g}_{\varepsilon^{(j)}(t_j'),\delta(t_j')}} \, d\text{\rm vol}_{\bar{g}_{\varepsilon^{(j)}(t_j'),\delta(t_j')}} = 0,\] 
where $t_j' := t_j+t$. Taking the limit as $j \to \infty$, we obtain 
\[\int_{[-\frac{1}{2},\frac{1}{2}]^4 \setminus \{0\}} \langle \hat{h}(t),g_{\text{\rm eucl}} \rangle_{g_{\text{\rm eucl}}} \, d\text{\rm vol}_{g_{\text{\rm eucl}}}\] 
by the dominated convergence theorem. 

Finally, since $\hat{h}$ is invariant under the maps $(x_1,x_2,x_3,x_4) \mapsto (x_2,-x_1,x_3,x_4)$, $(x_1,x_2,x_3,x_4) \mapsto (x_1,x_2,x_4,-x_3)$, and $(x_1,x_2,x_3,x_4) \mapsto (x_3,x_4,x_1,x_2)$, we conclude that $\int_{[-\frac{1}{2},\frac{1}{2}]^4 \setminus \{0\}} \hat{h}_{ii}(t) \, d\text{\rm vol}_{g_{\text{\rm eucl}}} = 0$ for all $i \in \{1,2,3,4\}$. This contradicts Proposition \ref{liouville.thm.3}. The proof of Theorem \ref{weighted.C0.estimate} is now complete.
\end{proof}

We next define suitable weighted H\"older spaces. 
To fix notation, we denote by $d_t(x,x')$ the Riemannian distance with respect to the metric $\bar{g}_{(-32\omega t)^{-\frac{1}{4}},(-t)^{-\frac{1}{400}}}$ from $x$ to $x'$. 
Moreover, we denote by $P_{x,x'}^t$ the parallel transport along a minimizing geodesic from $x$ to $x'$ with respect to the metric $\bar{g}_{(-32\omega t)^{-\frac{1}{4}},(-t)^{-\frac{1}{400}}}$. 

\begin{definition} 
Given real numbers $\alpha \in (0,1)$, $\gamma>0$, $\sigma>0$, and $\Lambda>0$, we define $\|h\|_{X_{\gamma,\sigma,\Lambda}^{0,\alpha}}$ to be the supremum of the quantity 
\begin{align*} 
&(-t)^\gamma \, ((-t)^{-\frac{1}{4}}+r)^\sigma \, |h(x,t)| \\ 
&+ (-t)^\gamma \, ((-t)^{-\frac{1}{4}}+r)^{\sigma+2\alpha} \, (d_t(x,x')^2+|t-t'|)^{-\alpha} \, |P_{x,x'}^t h(x,t)-h(x',t')|. 
\end{align*} 
Here, the supremum is taken over all numbers $r \in (0,10)$, all times $t,t' \in (-\infty,-\Lambda]$ satisfying $|t-t'| \leq ((-t)^{-\frac{1}{4}} + r)^2$, and all points $x,x' \in [-\frac{1}{2},\frac{1}{2}]^4$ satisfying $|x|,|x'| \in [\frac{1}{2} \, r,(-t)^{-\frac{1}{4}}+r]$. Moreover, the norm of $h(x,t)$ is taken with respect to the metric $\bar{g}_{(-32\omega t)^{-\frac{1}{4}},(-t)^{-\frac{1}{400}}}$. We next define 
\[\|h\|_{X_{\gamma,\sigma,\Lambda}^{l,\alpha}} := \|h\|_{X_{\gamma,\sigma,\Lambda}^{0,\alpha}} + \|Dh\|_{X_{\gamma,\sigma+1,\Lambda}^{0,\alpha}} + \|D^2 h\|_{X_{\gamma,\sigma+2,\Lambda}^{0,\alpha}} + \|\frac{\partial}{\partial t} h\|_{X_{\gamma,\sigma+2,\Lambda}^{0,\alpha}},\] 
where $D$ denotes the Riemannian connection with respect to $\bar{g}_{(-32\omega t)^{-\frac{1}{4}},(-t)^{-\frac{1}{400}}}$. 
\end{definition}
Note that the space decay rates in the above definition are appropriately adjusted for the derivatives so that they accommodate for the scaling involved. 
This allows us to combine Proposition \ref{weighted.C0.estimate} with standard interior estimates for parabolic equations to draw the following conclusion:

\begin{corollary}
\label{a.priori.estimate.1}
Given real numbers $\gamma>0$ and $\sigma \in (0,2)$, we can find real numbers $\Lambda>0$ and $C>0$ with the following property. Suppose that $\varepsilon(t)$ is a function which is defined on the interval $(-\infty,-\Lambda]$ and satisfies Assumption \ref{assumption.on.epsilon}, and let $\delta(t) = (-t)^{-\frac{1}{400}}$. Let $h \in X_{\gamma,\sigma,\Lambda}^{1,\alpha}$ be a solution of the inhomogeneous heat equation $\frac{\partial}{\partial t} h(t) = \Delta_{L,\bar{g}_{\varepsilon(t),\delta(t)}} h(t) + \psi(t)$ which is defined on $M \times (-\infty,-\Lambda]$. Moreover, we assume that $h$ is invariant under the group $\mathscr{G}$ and satisfies the orthogonality conditions 
\[\int_{[-\frac{1}{2},\frac{1}{2}]^4 \setminus \{0\}} \langle h(t),\bar{o}_{1,\varepsilon(t),\delta(t)} \rangle_{\bar{g}_{\varepsilon(t),\delta(t)}} \, d\text{\rm vol}_{\bar{g}_{\varepsilon(t),\delta(t)}} = 0\] 
and 
\[\int_{[-\frac{1}{2},\frac{1}{2}]^4 \setminus \{0\}} \langle h(t),\bar{g}_{\varepsilon(t),\delta(t)} \rangle_{\bar{g}_{\varepsilon(t),\delta(t)}} \, d\text{\rm vol}_{\bar{g}_{\varepsilon(t),\delta(t)}} = 0\] 
for all $t \in (-\infty,-\Lambda]$. Then 
\[\|h\|_{X_{\gamma,\sigma,\Lambda}^{1,\alpha}} \leq C \, \|\psi\|_{X_{\gamma,\sigma+2,\Lambda}^{0,\alpha}}.\] 
\end{corollary}
The following result is the main result of this section:

\begin{proposition}
\label{a.priori.estimate.2}
Given real numbers $\gamma>0$ and $\sigma \in (0,2)$, we can find real numbers $\Lambda>0$ and $C>0$ with the following property. Suppose that $\varepsilon(t)$ is a function which is defined on the interval $(-\infty,-\Lambda]$ and satisfies Assumption \ref{assumption.on.epsilon}, and let $\delta(t) = (-t)^{-\frac{1}{400}}$. Let $h \in X_{\gamma,\sigma,\Lambda}^{1,\alpha}$ be a solution of the inhomogeneous heat equation 
\[\frac{\partial}{\partial t} h(t) = \Delta_{L,\bar{g}_{\varepsilon(t),\delta(t)}} h(t) + \psi(t) + \lambda(t) \, \bar{o}_{1,\varepsilon(t),\delta(t)} + \nu(t) \, \bar{g}_{\varepsilon(t),\delta(t)}\] 
which is defined on $M \times (-\infty,-\Lambda]$. Moreover, we assume that $h$ is invariant under the group $\mathscr{G}$ and satisfies the orthogonality conditions 
\[\int_{[-\frac{1}{2},\frac{1}{2}]^4 \setminus \{0\}} \langle h(t),\bar{o}_{1,\varepsilon(t),\delta(t)} \rangle_{\bar{g}_{\varepsilon(t),\delta(t)}} \, d\text{\rm vol}_{\bar{g}_{\varepsilon(t),\delta(t)}} = 0\] 
and 
\[\int_{[-\frac{1}{2},\frac{1}{2}]^4 \setminus \{0\}} \langle h(t),\bar{g}_{\varepsilon(t),\delta(t)} \rangle_{\bar{g}_{\varepsilon(t),\delta(t)}} \, d\text{\rm vol}_{\bar{g}_{\varepsilon(t),\delta(t)}} = 0\] 
for all $t \in (-\infty,-\Lambda]$. Then 
\[\|h\|_{X_{\gamma,\sigma,\Lambda}^{1,\alpha}} \leq C \, \|\psi\|_{X_{\gamma+\alpha,\sigma+2,\Lambda}^{0,\alpha}}.\] 
Moreover, the function 
\[E(t) := \int_{[-\frac{1}{2},\frac{1}{2}]^4 \setminus \{0\}} \langle \psi(t) + \lambda(t) \, \bar{o}_{1,\varepsilon(t),\delta(t)},\bar{o}_{1,\varepsilon(t),\delta(t)} \rangle_{\bar{g}_{\varepsilon(t),\delta(t)}} \, d\text{\rm vol}_{\bar{g}_{\varepsilon(t),\delta(t)}}\] 
satisfies 
\[|E(t)| \leq C \, (-t)^{-1-\gamma-\alpha} \, \delta(t)^{2-\sigma} \, \|\psi\|_{X_{\gamma+\alpha,\sigma+2,\Lambda}^{0,\alpha}}\] 
and 
\[|t-t'|^{-\alpha} \, |E(t)-E(t')| \leq C \, (-t)^{-1-\gamma-\frac{\alpha}{2}} \, \delta(t)^{2-\sigma} \, \|\psi\|_{X_{\gamma+\alpha,\sigma+2,\Lambda}^{0,\alpha}}\]
for $0 < |t-t'| \leq (-t)^{-\frac{1}{2}}$. Finally, the function 
\[F(t) := \int_{[-\frac{1}{2},\frac{1}{2}]^4 \setminus \{0\}} \langle \psi(t) + \nu(t) \, \bar{g}_{\varepsilon(t),\delta(t)},\bar{g}_{\varepsilon(t),\delta(t)} \rangle_{\bar{g}_{\varepsilon(t),\delta(t)}} \, d\text{\rm vol}_{\bar{g}_{\varepsilon(t),\delta(t)}}\] 
satisfies 
\[|F(t)| \leq C \, (-t)^{-2-\gamma+\frac{\sigma}{4}} \, \|\psi\|_{X_{\gamma+\alpha,\sigma+2,\Lambda}^{0,\alpha}}\] 
and 
\[|t-t'|^{-\alpha} \, |F(t)-F(t')| \leq C \, (-t)^{-2-\gamma+\frac{\sigma}{4}+\frac{\alpha}{2}} \, \|\psi\|_{X_{\gamma+\alpha,\sigma+2,\Lambda}^{0,\alpha}}\] 
for $0 < |t-t'| \leq (-t)^{-\frac{1}{2}}$.
\end{proposition}

\begin{proof} 
Let us define 
\begin{align*} 
E_0(t) &= \int_{[-\frac{1}{2},\frac{1}{2}]^4 \setminus \{0\}} \langle \frac{\partial}{\partial t} h(t),\bar{o}_{1,\varepsilon(t),\delta(t)} \rangle_{\bar{g}_{\varepsilon(t),\delta(t)}} \, d\text{\rm vol}_{\bar{g}_{\varepsilon(t),\delta(t)}}, \\ 
E_1(t) &= \int_{[-\frac{1}{2},\frac{1}{2}]^4 \setminus \{0\}} \langle h(t),\Delta_{L,\bar{g}_{\varepsilon(t),\delta(t)}} \bar{o}_{1,\varepsilon(t),\delta(t)} \rangle_{\bar{g}_{\varepsilon(t),\delta(t)}} \, d\text{\rm vol}_{\bar{g}_{\varepsilon(t),\delta(t)}}, \\ 
E_2(t) &= \int_{[-\frac{1}{2},\frac{1}{2}]^4 \setminus \{0\}} \langle \psi,\bar{o}_{1,\varepsilon(t),\delta(t)} \rangle_{\bar{g}_{\varepsilon(t),\delta(t)}} \, d\text{\rm vol}_{\bar{g}_{\varepsilon(t),\delta(t)}}. 
\end{align*} 
Differentiating the identity 
\[\int_{[-\frac{1}{2},\frac{1}{2}]^4 \setminus \{0\}} \langle h(t),\bar{o}_{1,\varepsilon(t),\delta(t)} \rangle_{\bar{g}_{\varepsilon(t),\delta(t)}} \, d\text{\rm vol}_{\bar{g}_{\varepsilon(t),\delta(t)}} = 0\] 
with respect to $t$ gives 
\[|E_0(t)| \leq C \, (-t)^{-2-\gamma+\frac{\sigma}{4}} \, \|h\|_{X_{\gamma,\sigma,\Lambda}^{0,\alpha}}\] 
and 
\[|t-t'|^{-\alpha} \, |E_0(t)-E_0(t')| \leq C \, (-t)^{-2-\gamma+\frac{\sigma}{4}+\frac{\alpha}{2}} \, \|h\|_{X_{\gamma,\sigma,\Lambda}^{0,\alpha}}\] 
for $0 < |t-t'| \leq (-t)^{-\frac{1}{2}}$. Moreover, in view of Proposition \ref{obstruction.tensor}, the error term $E_1(t)$ satisfies 
\[|E_1(t)| \leq C \, (-t)^{-1-\gamma} \, \delta(t)^{2-\sigma} \, \|h\|_{X_{\gamma,\sigma,\Lambda}^{0,\alpha}}\] 
and 
\[|t-t'|^{-\alpha} \, |E_1(t)-E_1(t')| \leq C \, (-t)^{-1-\gamma+\frac{\alpha}{2}} \, \delta(t)^{2-\sigma} \, \|h\|_{X_{\gamma,\sigma,\Lambda}^{0,\alpha}}\]
for $0 < |t-t'| \leq (-t)^{-\frac{1}{2}}$. Finally, we have 
\[|E_2(t)| \leq C \, (-t)^{-\frac{1}{2}-\gamma+\frac{\sigma}{4}-\alpha} \, \|\psi\|_{X_{\gamma+\alpha,\sigma+2,\Lambda}^{0,\alpha}}\] 
and 
\[|t-t'|^{-\alpha} \, |E_2(t)-E_2(t')| \leq C \, (-t)^{-\frac{1}{2}-\gamma+\frac{\sigma}{4}-\frac{\alpha}{2}} \, \|\psi\|_{X_{\gamma+\alpha,\sigma+2,\Lambda}^{0,\alpha}}\] 
for $0 < |t-t'| \leq (-t)^{-\frac{1}{2}}$. Using the identity 
\[\lambda(t) \int_{[-\frac{1}{2},\frac{1}{2}]^4 \setminus \{0\}} |\bar{o}_{1,\varepsilon(t),\delta(t)}|_{\bar{g}_{\varepsilon(t),\delta(t)}}^2 \, d\text{\rm vol}_{\bar{g}_{\varepsilon(t),\delta(t)}} = E_0(t) - E_1(t) - E_2(t),\] 
we conclude that 
$$
|\lambda(t)| 
\leq C \, (-t)^{-\gamma} \, \delta(t)^{2-\sigma} \, \|h\|_{X_{\gamma,\sigma,\Lambda}^{0,\alpha}} 
+ C \, (-t)^{\frac{1}{2}-\gamma+\frac{\sigma}{4}-\alpha} \, \|\psi\|_{X_{\gamma+\alpha,\sigma+2,\Lambda}^{0,\alpha}} 
$$
and 
\begin{align*} 
|t-t'|^{-\alpha} \, |\lambda(t)-\lambda(t')| 
&\leq C \, (-t)^{-\gamma+\frac{\alpha}{2}} \, \delta(t)^{2-\sigma} \, \|h\|_{X_{\gamma,\sigma,\Lambda}^{0,\alpha}} \\ 
&+ C \, (-t)^{\frac{1}{2}-\gamma+\frac{\sigma}{4}-\frac{\alpha}{2}} \, \|\psi\|_{X_{\gamma+\alpha,\sigma+2,\Lambda}^{0,\alpha}} 
\end{align*}
for $0 < |t-t'| \leq (-t)^{-\frac{1}{2}}$. From this, we deduce that 
\[\|\lambda(\cdot) \, \bar{o}_{1,\varepsilon(\cdot),\delta(\cdot)}\|_{X_{\gamma,\sigma+2,\Lambda}^{0,\alpha}} \leq o(1) \, \|h\|_{X_{\gamma,\sigma,\Lambda}^{0,\alpha}} + C \, \|\psi\|_{X_{\gamma+\alpha,\sigma+2,\Lambda}^{0,\alpha}},\] 
where $o(1)$ represents a term that goes to $0$ as $\Lambda \to \infty$. 

In the next step, we define 
\begin{align*} 
F_0(t) &= \int_{[-\frac{1}{2},\frac{1}{2}]^4 \setminus \{0\}} \langle \frac{\partial}{\partial t} h(t),\bar{g}_{\varepsilon(t),\delta(t)} \rangle_{\bar{g}_{\varepsilon(t),\delta(t)}} \, d\text{\rm vol}_{\bar{g}_{\varepsilon(t),\delta(t)}}, \\ 
F_1(t) &= \int_{[-\frac{1}{2},\frac{1}{2}]^4 \setminus \{0\}} \langle \psi,\bar{g}_{\varepsilon(t),\delta(t)} \rangle_{\bar{g}_{\varepsilon(t),\delta(t)}}  \, d\text{\rm vol}_{\bar{g}_{\varepsilon(t),\delta(t)}}. 
\end{align*} 
Differentiating the relation 
\[\int_{[-\frac{1}{2},\frac{1}{2}]^4 \setminus \{0\}} \langle h(t),\bar{g}_{\varepsilon(t),\delta(t)} \rangle_{\bar{g}_{\varepsilon(t),\delta(t)}} \, d\text{\rm vol}_{\bar{g}_{\varepsilon(t),\delta(t)}} = 0\] 
with respect to $t$ yields 
\[|F_0(t)| \leq C \, (-t)^{-2-\gamma+\frac{\sigma}{4}} \, \|h\|_{X_{\gamma,\sigma,\Lambda}^{0,\alpha}}\] 
and 
\[|t-t'|^{-\alpha} \, |F_0(t)-F_0(t')| \leq C \, (-t)^{-2-\gamma+\frac{\sigma}{4}+\frac{\alpha}{2}} \, \|h\|_{X_{\gamma,\sigma,\Lambda}^{0,\alpha}}\] 
for $0 < |t-t'| \leq (-t)^{-\frac{1}{2}}$. We next observe that 
\[|F_1(t)| \leq C \, (-t)^{-\gamma-\alpha} \, \|\psi\|_{X_{\gamma+\alpha,\sigma+2,\Lambda}^{0,\alpha}}\] 
and 
\[|t-t'|^{-\alpha} \, |F_1(t)-F_1(t')| \leq C \, (-t)^{-\gamma-\frac{\alpha}{2}} \, \|\psi\|_{X_{\gamma+\alpha,\sigma+2,\Lambda}^{0,\alpha}}\] 
for $0 < |t-t'| \leq (-t)^{-\frac{1}{2}}$. Using the identity 
\[\nu(t) \int_{[-\frac{1}{2},\frac{1}{2}]^4 \setminus \{0\}} 4 \, d\text{\rm vol}_{\bar{g}_{\varepsilon(t),\delta(t)}} = F_0(t) - F_1(t),\] 
we conclude that 
$$
|\nu(t)| 
\leq C \, (-t)^{-2-\gamma+\frac{\sigma}{4}} \, \|h\|_{X_{\gamma,\sigma,\Lambda}^{0,\alpha}} 
+ C \, (-t)^{-\gamma-\alpha} \, \|\psi\|_{X_{\gamma+\alpha,\sigma+2,\Lambda}^{0,\alpha}} 
$$
and 
$$
|t-t'|^{-\alpha} \, |\nu(t)-\nu(t')| 
\leq C \, (-t)^{-2-\gamma+\frac{\sigma}{4}+\frac{\alpha}{2}} \, \|h\|_{X_{\gamma,\sigma,\Lambda}^{0,\alpha}} 
+ C \, (-t)^{-\gamma-\frac{\alpha}{2}} \, \|\psi\|_{X_{\gamma+\alpha,\sigma+2,\Lambda}^{0,\alpha}} 
$$
for $0 < |t-t'| \leq (-t)^{-\frac{1}{2}}$. This gives 
\[\|\nu(\cdot) \, \bar{g}_{\varepsilon(\cdot),\delta(\cdot)}\|_{X_{\gamma,\sigma+2,\Lambda}^{0,\alpha}} \leq o(1) \, \|h\|_{X_{\gamma,\sigma,\Lambda}^{0,\alpha}} + C \, \|\psi\|_{X_{\gamma+\alpha,\sigma+2,\Lambda}^{0,\alpha}},\] 
where again $o(1)$ represents a term that converges to $0$ as $\Lambda \to \infty$. 

After these preparations, we can now complete the proof. Using Corollary \ref{a.priori.estimate.1}, we obtain  
\[\|h\|_{X_{\gamma,\sigma,\Lambda}^{1,\alpha}} \leq C \, \|\psi\|_{X_{\gamma,\sigma+2,\Lambda}^{0,\alpha}} + C \, \|\lambda(\cdot) \, \bar{o}_{1,\varepsilon(\cdot),\delta(\cdot)}\|_{X_{\gamma,\sigma+2,\Lambda}^{0,\alpha}} + C \, \|\nu(\cdot) \, \bar{g}_{\varepsilon(\cdot),\delta(\cdot)}\|_{X_{\gamma,\sigma+2,\Lambda}^{0,\alpha}},\] 
hence
\[\|h\|_{X_{\gamma,\sigma,\Lambda}^{1,\alpha}} \leq o(1) \, \|h\|_{X_{\gamma,\sigma,\Lambda}^{0,\alpha}} + C \, \|\psi\|_{X_{\gamma+\alpha,\sigma+2,\Lambda}^{0,\alpha}}\] 
if $\Lambda$ is sufficiently large. From this, the first statement follows. Finally, the estimates for $E(t)$ and $F(t)$ follow from the fact that $E(t) = E_0(t)-E_1(t)$ and $F(t) = F_0(t)$.
\end{proof}

\begin{corollary} 
\label{existence.uniqueness}
Let $\gamma>0$ and $\sigma \in (0,2)$ be arbitrary, and let $\Lambda>0$ be chosen as in Proposition \ref{a.priori.estimate.2}. Then, given any tensor $\psi \in X_{\gamma+\alpha,\sigma+2,\Lambda}^{0,\alpha}$, there exists a unique tensor $h \in X_{\gamma,\sigma,\Lambda}^{1,\alpha}$ and scalar functions $\lambda: (-\infty,-\Lambda] \to \mathbb{R}$ and $\nu: (-\infty,-\Lambda] \to \mathbb{R}$ such that 
\[\frac{\partial}{\partial t} h(t) = \Delta_{L,\bar{g}_{\varepsilon(t),\delta(t)}} h(t) + \psi(t) + \lambda(t) \, \bar{o}_{1,\varepsilon(t),\delta(t)} + \nu(t) \, \bar{g}_{\varepsilon(t),\delta(t)}.\] 
Moreover, $h$ is invariant under the group $\mathscr{G}$, and satisfies the orthogonality conditions 
\[\int_{[-\frac{1}{2},\frac{1}{2}]^4 \setminus \{0\}} \langle h(t),\bar{o}_{1,\varepsilon(t),\delta(t)} \rangle_{\bar{g}_{\varepsilon(t),\delta(t)}} \, d\text{\rm vol}_{\bar{g}_{\varepsilon(t),\delta(t)}} = 0\] 
and 
\[\int_{[-\frac{1}{2},\frac{1}{2}]^4 \setminus \{0\}} \langle h(t),\bar{g}_{\varepsilon(t),\delta(t)} \rangle_{\bar{g}_{\varepsilon(t),\delta(t)}} \, d\text{\rm vol}_{\bar{g}_{\varepsilon(t),\delta(t)}} = 0\] 
for all $t \in (-\infty,-\Lambda]$.
\end{corollary}

\begin{proof} 
The uniqueness statement follows immediately from Proposition \ref{a.priori.estimate.2}. Hence, it remains to prove the existence statement. If $\psi$ is compactly supported, standard results on linear parabolic equations imply that there exists a compactly supported tensor $h$ and scalar functions $\lambda(\cdot)$ and $\nu(\cdot)$ with the required properties. To prove the assertion in general, we approximate $\psi$ by compactly supported tensors, and use the a priori estimate in Proposition \ref{a.priori.estimate.2} to pass to the limit.
\end{proof}

\section{Existence of a solution to the nonlinear problem}

Throughout this section, we fix positive real numbers $\alpha,\gamma,\sigma \in (0,1)$. We assume that $\alpha$ and $\sigma$ are very small, and $\gamma$ is close to $1$. Specifically, we may choose $\alpha \in (0,\frac{1}{10000})$, $\gamma \in (1-\frac{1}{10000},1)$, and $\sigma \in (0,\frac{1}{10000})$. In addition, we require that $\gamma+\alpha < 1$. Furthermore, $\Lambda$ will denote a positive real number which we will choose sufficiently large depending on $\alpha,\gamma,\sigma$. 

Given a metric $g$ and a symmetric $(0,2)$-tensor $k$ satisfying $|k|_g \leq \frac{1}{2}$, we define 
\[Q_g(k) = 2 \, \text{\rm Ric}_{g+k} - 2 \, \text{\rm Ric}_g + \Delta_{L,g} k - \mathscr{L}_Y (g+k),\] 
where $Y = \text{\rm div}_g k - \frac{1}{2} \, \nabla \text{\rm tr}_g k$. Note that $Q_g(k)$ can be expanded as 
\[Q_g(k) = k * \nabla^2 k + \nabla k * \nabla k + \text{\rm Rm}_g * k * k + \text{\rm higher order terms.}\] 

Recall the definition of $\bar{g}_{\varepsilon,\delta}$ from Section \ref{attaching}. Our goal is to perturb the family of metrics $\bar{g}_{(-32\omega t)^{-\frac{1}{4}},(-t)^{-\frac{1}{400}}}$ to an exact solution to the Ricci flow which is defined for $t \in (-\infty,-\Lambda]$. This problem comes down to finding a fixed point of a nonlinear mapping between Banach spaces. In the following, we describe this mapping in detail. Let $\mathscr{A}_{\gamma,\sigma,\Lambda}^\alpha$ denote the set of all triplets $(k,\eta(\cdot),\beta(\cdot))$ which satisfy the following conditions: 
\begin{itemize} 
\item The tensor $k$ is invariant under $\mathscr{G}$ and satisfies $\|k\|_{X_{\gamma,\sigma,\Lambda}^{1,\alpha}} \leq 1$. 
\item The function $\eta: (-\infty,-\Lambda] \to \mathbb{R}$ satisfies $|\eta(t)| \leq (-t)^{-\frac{1}{1000}}$, 
\\ and $|t-t'|^{-\alpha} \, |\eta(t)-\eta(t')| \leq (-t)^{-\frac{1}{1000}}$ for $0 < |t-t'| \leq (-t)^{-\frac{1}{2}}$. 
\item The function $\beta: (-\infty,-\Lambda] \to \mathbb{R}$ satisfies $|\beta(t)| \leq (-t)^{-1}$ 
\\ and $|t-t'|^{-\alpha} \, |\beta(t)-\beta(t')| \leq (-t)^{-1}$ for $0 < |t-t'| \leq (-t)^{-\frac{1}{2}}$. 
\end{itemize}
We can think of $\mathscr{A}_{\gamma,\sigma,\Lambda}^\alpha$ as the unit ball in a suitable Banach space. This Banach space will be denoted by $\mathscr{E}_{\gamma,\sigma,\Lambda}^\alpha$. 

Given a triplet $(k,\eta(\cdot),\beta(\cdot)) \in \mathscr{A}_{\gamma,\sigma,\Lambda}^\alpha$, we consider the family of metrics $\bar{g}_{\varepsilon(t),\delta(t)}$, where 
\[\varepsilon(t) = \bigg ( - 32\omega t + \int_t^{-\Lambda} \eta(s) \, ds \bigg )^{-\frac{1}{4}}\]  
and $\delta(t) = (-t)^{-\frac{1}{400}}$ for $t \in (-\infty,-\Lambda]$. It is straightforward to verify that the function $\varepsilon(t)$ satisfies Assumption \ref{assumption.on.epsilon} provided that $\Lambda$ is sufficiently large. By Corollary \ref{existence.uniqueness}, there exists a unique tensor $h \in X_{\gamma,\sigma,\Lambda}^{1,\alpha}$ and scalar functions $\lambda: (-\infty,-\Lambda] \to \mathbb{R}$ and $\nu: (-\infty,-\Lambda] \to \mathbb{R}$ such that 
\begin{gather*} 
\frac{\partial}{\partial t} h(t) 
= \Delta_{L,\bar{g}_{\varepsilon(t),\delta(t)}} h(t) + \psi(t) 
+ \lambda(t) \, \bar{o}_{1,\varepsilon(t),\delta(t)} + \nu(t) \, \bar{g}_{\varepsilon(t),\delta(t)}, 
\\
\text{ where } 
\psi(t) := -Q_{\bar{g}_{\varepsilon(t),\delta(t)}}(k(t)) + \beta(t) \, k(t) - \Big ( \frac{\partial}{\partial t} \bar{g}_{\varepsilon(t),\delta(t)} + 2 \, \text{\rm Ric}_{\bar{g}_{\varepsilon(t),\delta(t)}} \Big ),  
\end{gather*} 
and $h(t)$ satisfies the orthogonality conditions 
\[\int_{[-\frac{1}{2},\frac{1}{2}]^4 \setminus \{0\}} \langle h(t),\bar{o}_{1,\varepsilon(t),\delta(t)} \rangle_{\bar{g}_{\varepsilon(t),\delta(t)}} \, d\text{\rm vol}_{\bar{g}_{\varepsilon(t),\delta(t)}} = 0\] 
and 
\[\int_{[-\frac{1}{2},\frac{1}{2}]^4 \setminus \{0\}} \langle h(t),\bar{g}_{\varepsilon(t),\delta(t)} \rangle_{\bar{g}_{\varepsilon(t),\delta(t)}} \, d\text{\rm vol}_{\bar{g}_{\varepsilon(t),\delta(t)}} = 0\] 
for all $t \in (-\infty,-\Lambda]$. 

Having solved this linear PDE for $h$, $\lambda(\cdot)$, $\nu(\cdot)$, we next define a function $\xi: (-\infty,-\Lambda] \to \mathbb{R}$ by 
\[\xi(t) = \eta(t) - \pi^{-2} \, \varepsilon(t)^{-8} \, \lambda(t) \int_{[-\frac{1}{2},\frac{1}{2}]^4 \setminus \{0\}} |\bar{o}_{1,\varepsilon(t),\delta(t)}|_{\bar{g}_{\varepsilon(t),\delta(t)}}^2 \, d\text{\rm vol}_{\bar{g}_{\varepsilon(t),\delta(t)}}\] 
for $t \in (-\infty,-\Lambda]$.

\textbf{In the remainder of this section, we will analyze the map $\mathscr{J}$ which sends the triplet $(k,\eta(\cdot),\beta(\cdot)) \in \mathscr{A}_{\gamma,\sigma,\Lambda}^\alpha$ to the triplet $(h,\xi(\cdot),\nu(\cdot))$.} In order to prove the existence of a fixed point, we need to show that $\mathscr{J}$ maps the unit ball $\mathscr{A}_{\gamma,\sigma,\Lambda}^\alpha$ into itself. By combining Corollary \ref{projection.of.Ric.to.o_1}, Proposition \ref{projection.of.Ric.to.g}, and Proposition \ref{projection.of.time.derivative}, we obtain the following result, which relates the time derivative of $\varepsilon(t)$ to the orthogonal projection of $\psi$ to the approximate kernel. This plays the role of a balancing condition; it serves as the main motivation for the definition of $\varepsilon(t)$ above. 

\begin{proposition}
\label{estimate.for.psi}
Consider a triplet $(k,\eta(\cdot),\beta(\cdot)) \in \mathscr{A}_{\gamma,\sigma,\Lambda}^\alpha$, and let 
$\psi(t)$ and $h(t)$ be as in the discussion above. 
Then $\|\psi\|_{X_{\gamma+\alpha,\sigma+2,\Lambda}^{0,\alpha}} \leq o(1)$, where $o(1)$ represents a term that converges to $0$ as $\Lambda \to \infty$. Moreover, the function 
\begin{align*} 
G(t) 
&:= 4\pi^2 \, \varepsilon(t)^3 \, \frac{d\varepsilon}{dt}(t) - 32\pi^2 \, \omega \, \varepsilon(t)^8 \\ 
&+ \int_{[-\frac{1}{2},\frac{1}{2}]^4 \setminus \{0\}} \langle \psi(t),\bar{o}_{1,\varepsilon(t),\delta(t)} \rangle_{\bar{g}_{\varepsilon(t),\delta(t)}} \, d\text{\rm vol}_{\bar{g}_{\varepsilon(t),\delta(t)}} 
\end{align*}
satisfies $|G(t)| \leq C \, (-t)^{-2-\frac{1}{400}}$ and $|t-t'|^{-\alpha} \, |G(t)-G(t')| \leq C \, (-t)^{-2-\frac{1}{400}}$ for $0 < |t-t'| \leq (-t)^{-\frac{1}{2}}$. Finally, the function 
\[H(t) := \int_{[-\frac{1}{2},\frac{1}{2}]^4 \setminus \{0\}} \langle \psi(t),\bar{g}_{\varepsilon(t),\delta(t)} \rangle_{\bar{g}_{\varepsilon(t),\delta(t)}} \, d\text{\rm vol}_{\bar{g}_{\varepsilon(t),\delta(t)}}\] 
satisfies $|H(t)| \leq C \, (-t)^{-\frac{5}{4}}$ and $|t-t'|^{-\alpha} \, |H(t)-H(t')| \leq C \, (-t)^{-\frac{5}{4}}$ for $0 < |t-t'| \leq (-t)^{-\frac{1}{2}}$.
\end{proposition}

\begin{proof} 
The inequality $\|k\|_{X_{\gamma,\sigma,\Lambda}^{1,\alpha}} \leq 1$ implies $\|Q_g(k)\|_{X_{\gamma+\alpha,\sigma+2,\Lambda}^{0,\alpha}} \leq o(1)$ and $\|\beta \, k\|_{X_{\gamma+\alpha,\sigma+2,\Lambda}^{0,\alpha}} \leq o(1)$. Moreover, since $\gamma+\alpha < 1$, we have 
\[\big \| \frac{\partial}{\partial t} \bar{g}_{\varepsilon(t),\delta(t)} + 2 \, \text{\rm Ric}_{\bar{g}_{\varepsilon(t),\delta(t)}} \big \|_{X_{\gamma+\alpha,\sigma+2,\Lambda}^{0,\alpha}} \leq o(1).\] 
Putting these facts together, we conclude that $\|\psi\|_{X_{\gamma+\alpha,\sigma+2,\Lambda}^{0,\alpha}} \leq o(1)$.
 
We next estimate the function $G(t)$. Using the inequality $\|k\|_{X_{\gamma,\sigma,\Lambda}^{1,\alpha}} \leq 1$, we obtain $|Q_{\bar{g}_{\varepsilon(t),\delta(t)}}(k(t))| \leq C \, (-t)^{-2\gamma} \, ((-t)^{-\frac{1}{4}} + |x|)^{-2-2\sigma}$, hence 
\begin{align*} 
&\bigg | \int_{[-\frac{1}{2},\frac{1}{2}]^4 \setminus \{0\}} \langle Q_{\bar{g}_{\varepsilon(t),\delta(t)}}(k(t)),\bar{o}_{1,\varepsilon(t),\delta(t)} \rangle_{\bar{g}_{\varepsilon(t),\delta(t)}} \, d\text{\rm vol}_{\bar{g}_{\varepsilon(t),\delta(t)}} 
\bigg | \\ 
&\leq C \int_{[-\frac{1}{2},\frac{1}{2}]^4 \setminus \{0\}} (-t)^{-1-2\gamma} \, ((-t)^{-\frac{1}{4}}+|x|)^{-6-2\sigma} \, d\text{\rm vol}_{\bar{g}_{\varepsilon(t),\delta(t)}} \\ 
&\leq C \, (-t)^{-\frac{1}{2}-2\gamma+\frac{\sigma}{2}}. 
\end{align*}
Moreover, since $|k| \leq (-t)^{-\gamma} \, ((-t)^{-\frac{1}{4}}+|x|)^{-\sigma}$ and $|\beta(t)| \leq (-t)^{-1}$, we obtain  
\begin{align*} 
&\bigg | \int_{[-\frac{1}{2},\frac{1}{2}]^4 \setminus \{0\}} \langle \beta(t) \, k(t),\bar{o}_{1,\varepsilon(t),\delta(t)} \rangle_{\bar{g}_{\varepsilon(t),\delta(t)}} \, d\text{\rm vol}_{\bar{g}_{\varepsilon(t),\delta(t)}} \bigg | \\ 
&\leq C \int_{[-\frac{1}{2},\frac{1}{2}]^4 \setminus \{0\}} (-t)^{-2-\gamma} \, ((-t)^{-\frac{1}{4}}+|x|)^{-4-\sigma} \, d\text{\rm vol}_{\bar{g}_{\varepsilon(t),\delta(t)}} \\ 
&\leq C \, (-t)^{-2-\gamma+\frac{\sigma}{4}}. 
\end{align*}
Finally, using Corollary \ref{projection.of.Ric.to.o_1} and Proposition \ref{projection.of.time.derivative}, we obtain 
\begin{align*} 
&\bigg | 4\pi^2 \, \varepsilon(t)^3 \, \frac{d\varepsilon}{dt}(t) - 32\pi^2 \, \omega \, \varepsilon(t)^8 \\ 
&- \int_{[-\frac{1}{2},\frac{1}{2}]^4 \setminus \{0\}} \Big \langle \frac{\partial}{\partial t} \bar{g}_{\varepsilon(t),\delta(t)} + 2 \, \text{\rm Ric}_{\bar{g}_{\varepsilon(t),\delta(t)}},\bar{o}_{1,\varepsilon(t),\delta(t)} \Big \rangle_{\bar{g}_{\varepsilon(t),\delta(t)}} \, d\text{\rm vol}_{\bar{g}_{\varepsilon(t),\delta(t)}} \bigg | \\ 
&\leq C \, \varepsilon(t)^8 \, \delta(t)^2.
\end{align*}
Putting these facts together, we obtain $|G(t)| \leq C \, (-t)^{-2-\frac{1}{400}}$, as claimed. A similar argument gives $|t-t'|^{-\alpha} \, |G(t)-G(t')| \leq C \, (-t)^{-2-\frac{1}{400}}$ for $0 < |t-t'| \leq (-t)^{-\frac{1}{2}}$.

It remains to estimate the function $H(t)$. Using the inequality $\|k\|_{X_{\gamma,\sigma,\Lambda}^{1,\alpha}} \leq	1$, we obtain  
\begin{align*} 
&\bigg | \int_{[-\frac{1}{2},\frac{1}{2}]^4 \setminus \{0\}} \langle Q_{\bar{g}_{\varepsilon(t),\delta(t)}}(k(t)),\bar{g}_{\varepsilon(t),\delta(t)} \rangle_{\bar{g}_{\varepsilon(t),\delta(t)}} \, d\text{\rm vol}_{\bar{g}_{\varepsilon(t),\delta(t)}} 
\bigg | \\ 
&\leq C \int_{[-\frac{1}{2},\frac{1}{2}]^4 \setminus \{0\}} (-t)^{-2\gamma} \, ((-t)^{-\frac{1}{4}}+|x|)^{-2-2\sigma} \, d\text{\rm vol}_{\bar{g}_{\varepsilon(t),\delta(t)}} \\ 
&\leq C \, (-t)^{-2\gamma} 
\end{align*}
and 
\begin{align*} 
&\bigg | \int_{[-\frac{1}{2},\frac{1}{2}]^4 \setminus \{0\}} \langle \beta(t) \, k(t),\bar{g}_{\varepsilon(t),\delta(t)} \rangle_{\bar{g}_{\varepsilon(t),\delta(t)}} \, d\text{\rm vol}_{\bar{g}_{\varepsilon(t),\delta(t)}} 
\bigg | \\ 
&\leq C \int_{[-\frac{1}{2},\frac{1}{2}]^4 \setminus \{0\}} (-t)^{-1-\gamma} \, ((-t)^{-\frac{1}{4}}+|x|)^{-\sigma} \, d\text{\rm vol}_{\bar{g}_{\varepsilon(t),\delta(t)}} \\ 
&\leq C \, (-t)^{-1-\gamma}.
\end{align*}
Moreover, it follows from Proposition \ref{projection.of.Ric.to.g} and Proposition \ref{projection.of.time.derivative} that 
\begin{align*} 
&\bigg | \int_{[-\frac{1}{2},\frac{1}{2}]^4 \setminus \{0\}} \Big \langle \frac{\partial}{\partial t} \bar{g}_{\varepsilon(t),\delta(t)} + 2 \, \text{\rm Ric}_{\bar{g}_{\varepsilon(t),\delta(t)}},\bar{g}_{\varepsilon(t),\delta(t)} \Big \rangle_{\bar{g}_{\varepsilon(t),\delta(t)}} \, d\text{\rm vol}_{\bar{g}_{\varepsilon(t),\delta(t)}} \bigg | \\ 
&\leq C \, \varepsilon(t)^8 \, \delta(t)^{-6}.
\end{align*}
Putting these facts together, we obtain $|H(t)| \leq C \, (-t)^{-\frac{5}{4}}$. A similar argument gives $|t-t'|^{-\alpha} \, |H(t)-H(t')| \leq C \, (-t)^{-\frac{5}{4}}$ for $0 < |t-t'| \leq (-t)^{-\frac{1}{2}}$. This completes the proof.
\end{proof}

Combining Proposition \ref{estimate.for.psi} with Proposition \ref{a.priori.estimate.2}, we can draw the following conclusion: 

\begin{corollary} 
\label{image.of.J}
Consider a triplet $(k,\eta(\cdot),\beta(\cdot)) \in \mathscr{A}_{\gamma,\sigma,\Lambda}^\alpha$. Then: 
\begin{itemize} 
\item The tensor $h$ satisfies $\|h\|_{X_{\gamma,\sigma,\Lambda}^{1,\alpha}} \leq o(1)$. 
\item The function $\xi(\cdot)$ satisfies $|\xi(t)| \leq C \, (-t)^{-\frac{1}{400}}$ and $|t-t'|^{-\alpha} \, |\xi(t)-\xi(t')| \leq C \, (-t)^{-\frac{1}{400}}$ for $0 < |t-t'| \leq (-t)^{-\frac{1}{2}}$.
\item The function $\nu(\cdot)$ satisfies $|\nu(t)| \leq C \, (-t)^{-\frac{5}{4}}$ and $|t-t'|^{-\alpha} \, |\nu(t)-\nu(t')| \leq C \, (-t)^{-\frac{5}{4}}$ for $0 < |t-t'| \leq (-t)^{-\frac{1}{2}}$. 
\end{itemize}
Here, $C$ is a positive constant which does not depend on $\Lambda$, and $o(1)$ represents a quantity which converges to $0$ as $\Lambda \to \infty$.
\end{corollary}

\begin{proof}
For $k$ and $h$ as in the discussion preceding Proposition \ref{estimate.for.psi} 
we have by 
Proposition \ref{a.priori.estimate.2} that $\|h\|_{X_{\gamma,\sigma,\Lambda}^{1,\alpha}} \leq C \, \|\psi\|_{X_{\gamma+\alpha,\sigma+2,\Lambda}^{0,\alpha}}$.  
This implies by Proposition \ref{estimate.for.psi} that $\|h\|_{X_{\gamma,\sigma,\Lambda}^{1,\alpha}} \leq o(1)$ which proves the first statement. 

We next define  
\begin{align*} 
E(t) &= \int_{[-\frac{1}{2},\frac{1}{2}]^4 \setminus \{0\}} \langle \psi(t) + \lambda(t) \, \bar{o}_{1,\varepsilon(t),\delta(t)},\bar{o}_{1,\varepsilon(t),\delta(t)} \rangle_{\bar{g}_{\varepsilon(t),\delta(t)}} \, d\text{\rm vol}_{\bar{g}_{\varepsilon(t),\delta(t)}}, \\ 
F(t) &= \int_{[-\frac{1}{2},\frac{1}{2}]^4 \setminus \{0\}} \langle \psi(t) + \nu(t) \, \bar{g}_{\varepsilon(t),\delta(t)},\bar{g}_{\varepsilon(t),\delta(t)} \rangle_{\bar{g}_{\varepsilon(t),\delta(t)}} \, d\text{\rm vol}_{\bar{g}_{\varepsilon(t),\delta(t)}}, \\ 
G(t) &= 4\pi^2 \, \varepsilon(t)^3 \, \frac{d\varepsilon}{dt}(t) - 32\pi^2 \, \omega \, \varepsilon(t)^8 \\ 
&+ \int_{[-\frac{1}{2},\frac{1}{2}]^4 \setminus \{0\}} \langle \psi(t),\bar{o}_{1,\varepsilon(t),\delta(t)} \rangle_{\bar{g}_{\varepsilon(t),\delta(t)}} \, d\text{\rm vol}_{\bar{g}_{\varepsilon(t),\delta(t)}}, \\ 
H(t) &= \int_{[-\frac{1}{2},\frac{1}{2}]^4 \setminus \{0\}} \langle \psi(t),\bar{g}_{\varepsilon(t),\delta(t)} \rangle_{\bar{g}_{\varepsilon(t),\delta(t)}} \, d\text{\rm vol}_{\bar{g}_{\varepsilon(t),\delta(t)}}. 
\end{align*} 
By definition of $\xi(\cdot)$, we have 
\begin{align*} 
\xi(t) 
&= \eta(t) - \pi^{-2} \, \varepsilon(t)^{-8} \, \lambda(t) \int_{[-\frac{1}{2},\frac{1}{2}]^4 \setminus \{0\}} |\bar{o}_{1,\varepsilon(t),\delta(t)}|_{\bar{g}_{\varepsilon(t),\delta(t)}}^2 \, d\text{\rm vol}_{\bar{g}_{\varepsilon(t),\delta(t)}} \\ 
&= 4 \, \varepsilon(t)^{-5} \, \frac{d\varepsilon}{dt}(t) - 32\omega \\ 
&- \pi^{-2} \, \varepsilon(t)^{-8} \, \lambda(t) \int_{[-\frac{1}{2},\frac{1}{2}]^4 \setminus \{0\}} |\bar{o}_{1,\varepsilon(t),\delta(t)}|_{\bar{g}_{\varepsilon(t),\delta(t)}}^2 \, d\text{\rm vol}_{\bar{g}_{\varepsilon(t),\delta(t)}} \\ 
&= \pi^{-2} \, \varepsilon(t)^{-8} \, (G(t) - E(t)). 
\end{align*} 
Using Proposition \ref{a.priori.estimate.2} and Proposition \ref{estimate.for.psi}, we obtain $|\xi(t)| \leq C \, (-t)^{-\frac{1}{400}}$ and $|t-t'|^{-\alpha} \, |\xi(t)-\xi(t')| \leq C \, (-t)^{-\frac{1}{400}}$ for $0 < |t-t'| \leq (-t)^{-\frac{1}{2}}$. This proves the second statement.

Finally, we have 
\[\nu(t) \int_{[-\frac{1}{2},\frac{1}{2}]^4 \setminus \{0\}} 4 \, d\text{\rm vol}_{\bar{g}_{\varepsilon(t),\delta(t)}} = F(t) - H(t).\]
Hence, it follows from Proposition \ref{a.priori.estimate.2} and Proposition \ref{estimate.for.psi} that $|\nu(t)| \leq C \, (-t)^{-\frac{5}{4}}$ and $|t-t'|^{-\alpha} \, |\nu(t)-\nu(t')| \leq C \, (-t)^{-\frac{5}{4}}$ for $0 < |t-t'| \leq (-t)^{-\frac{1}{2}}$. From this, the third statement follows.
\end{proof}

After these preparations, we now prove that $\mathscr{J}$ has a fixed point. In addition to the original parameters $\gamma,\alpha$, we will consider another pair of parameters $\tilde{\gamma},\tilde{\alpha}$ such that $\tilde{\gamma} < \gamma$ and $\tilde{\alpha} < \alpha$.

\begin{proposition}
\label{properties.of.J}
The set $\mathscr{A}_{\gamma,\sigma,\Lambda}^\alpha$ is a compact subset of $\mathscr{E}_{\tilde{\gamma},\sigma,\Lambda}^{\tilde{\alpha}}$. Moreover, if $\Lambda$ is sufficiently large, then $\mathscr{J}$ maps the set $\mathscr{A}_{\gamma,\sigma,\Lambda}^\alpha$ into itself. Finally, the map $\mathscr{J}: \mathscr{A}_{\gamma,\sigma,\Lambda}^\alpha \to \mathscr{A}_{\gamma,\sigma,\Lambda}^\alpha$ is continuous with respect to the norm on $\mathscr{E}_{\tilde{\gamma},\sigma,\Lambda}^{\tilde{\alpha}}$. 
\end{proposition}

\begin{proof} 
The first statement is standard. The second statement follows directly from Corollary \ref{image.of.J}. In order to prove the third statement, we argue by contradiction. Suppose that there exists a sequence of triplets $u^{(j)} \in \mathscr{A}_{\gamma,\sigma,\Lambda}^\alpha$ and a triplet $u \in \mathscr{A}_{\gamma,\sigma,\Lambda}^\alpha$ such that $\limsup_{j \to \infty} \|u^{(j)}-u\|_{\mathscr{E}_{\tilde{\gamma},\sigma,\Lambda}^{\tilde{\alpha}}} = 0$ and $\liminf_{j \to \infty} \|\mathscr{J}(u^{(j)})-\mathscr{J}(u)\|_{\mathscr{E}_{\tilde{\gamma},\sigma,\Lambda}^{\tilde{\alpha}}} > 0$. Note that $\mathscr{A}_{\gamma,\sigma,\Lambda}^\alpha$ is a compact subset of $\mathscr{E}_{\tilde{\gamma},\sigma,\Lambda}^{\tilde{\alpha}}$ which contains the sequence $\mathscr{J}(u^{(j)})$. Hence, after passing to a subsequence, we can find an element $v \in \mathscr{A}_{\gamma,\sigma,\Lambda}^\alpha$ such that $\limsup_{j \to \infty} \|\mathscr{J}(u^{(j)})-v\|_{\mathscr{E}_{\tilde{\gamma},\sigma,\Lambda}^{\tilde{\alpha}}} = 0$. From the definition of $\mathscr{J}$, it is easy to see that $v = \mathscr{J}(u)$. This is a contradiction.
\end{proof}

\begin{corollary} 
\label{existence.of.fixed.point}
The map $\mathscr{J}: \mathscr{A}_{\gamma,\sigma,\Lambda}^\alpha \to \mathscr{A}_{\gamma,\sigma,\Lambda}^\alpha$ has a fixed point.
\end{corollary}

\begin{proof} 
This follows immediately from Proposition \ref{properties.of.J} and the Schauder fixed point theorem.
\end{proof}

\begin{proposition}
\label{ancient.solution}
Consider the map which sends a triplet $(k,\eta(\cdot),\beta(\cdot))$ to the triplet $(h,\xi(\cdot),\nu(\cdot))$. Every fixed point of this map corresponds to a solution of the Ricci flow on $M$ which is defined for $t \in (-\infty,-\Lambda]$.
\end{proposition}

\begin{proof} 
Let us consider a triplet $(k,\eta(\cdot),\beta(\cdot)) \in \mathscr{A}_{\gamma,\sigma,\Lambda}^\alpha$ with the property that $h=k$, $\xi(\cdot)=\eta(\cdot)$, and $\nu(\cdot)=\beta(\cdot)$. As above, we define 
\[\varepsilon(t) = \bigg ( - 32\omega t + \int_t^{-\Lambda} \eta(s) \, ds \bigg )^{-\frac{1}{4}}\] 
and $\delta(t) = (-t)^{-\frac{1}{400}}$ for $t \in (-\infty,-\Lambda]$. The relation $\xi(\cdot)=\eta(\cdot)$ directly implies $\lambda(\cdot)=0$. Moreover, since $h=k$ and $\nu(\cdot)=\beta(\cdot)$, it follows that $k$ is a solution of the equation 
\begin{align*} 
\frac{\partial}{\partial t} k(t) 
&= \Delta_{L,\bar{g}_{\varepsilon(t),\delta(t)}} k(t) - Q_{\bar{g}_{\varepsilon(t),\delta(t)}}(k(t)) + \beta(t) \, k(t) \\ 
&- \Big ( \frac{\partial}{\partial t} \bar{g}_{\varepsilon(t),\delta(t)} + 2 \, \text{\rm Ric}_{\bar{g}_{\varepsilon(t),\delta(t)}} \Big ) + \beta(t) \, \bar{g}_{\varepsilon(t),\delta(t)}.
\end{align*} 
Rearranging terms gives 
\begin{align*} 
&\frac{\partial}{\partial t} (\bar{g}_{\varepsilon(t),\delta(t)} + k(t)) \\ 
&= -2 \, \text{\rm Ric}_{\bar{g}_{\varepsilon(t),\delta(t)}} + \Delta_{L,\bar{g}_{\varepsilon(t),\delta(t)}} k(t) - Q_{\bar{g}_{\varepsilon(t),\delta(t)}}(k(t)) + \beta(t) \, (\bar{g}_{\varepsilon(t),\delta(t)} + k(t)) \\ 
&= -2 \, \text{\rm Ric}_{\bar{g}_{\varepsilon(t),\delta(t)}+k(t)} + \mathscr{L}_Y (\bar{g}_{\varepsilon(t),\delta(t)} + k(t)) + \beta(t) \, (\bar{g}_{\varepsilon(t),\delta(t)} + k(t)), 
\end{align*} 
where $Y = \text{\rm div}_{\bar{g}_{\varepsilon(t),\delta(t)}} k - \frac{1}{2} \nabla \text{\rm tr}_{\bar{g}_{\varepsilon(t),\delta(t)}} k$. Hence, if we put 
\[g(t) := e^{\int_t^{-\Lambda} \beta(s) \, ds} \, (\bar{g}_{\varepsilon(t),\delta(t)} + k(t)),\] 
then the metrics $g(t)$ satisfy 
\[\frac{\partial}{\partial t} g(t) = -2 \, \text{\rm Ric}_{g(t)} + \mathscr{L}_Y g(t).\] 
By pulling back the metrics $g(t)$ under the flow of diffeomorphisms generated by $Y(t)$, we obtain a solution to the Ricci flow which is defined for $t \in (-\infty,-\Lambda]$.
\end{proof}

\textit{Proof of Theorem \ref{main.thm}.} 
By Corollary \ref{existence.of.fixed.point}, the map $\mathscr{J}$ has a fixed point. By Proposition \ref{ancient.solution}, this corresponds to an ancient solution of the Ricci flow. \qed

\end{document}